\documentclass[12pt]{article}
\usepackage{amssymb}
\usepackage{amsxtra}
\usepackage{amstext}
\usepackage{latexsym}
\usepackage{amsmath}
\usepackage{dsfont} 

\newtheorem{definition}{Definition}[section]
\newtheorem{theorem}[definition]{Theorem}
\newtheorem{lemma}[definition]{Lemma}
\newtheorem{corollary}[definition]{Corollary}
\newtheorem{note}[definition]{Note}
\newtheorem{assumption}[definition]{Assumption}

\def\K{\mathbb K}
\def\Z{\mathbb Z}

\newcommand{\SL}{\mathfrak{sl}_2}
\newcommand{\uq}{U_q(\mathfrak{sl}_2)}
\newcommand{\Ab}{A^{*}}
\newcommand{\Vb}{V^{*}}
\newcommand{\Tb}{\theta^{*}}
\newcommand{\Eb}{E^{*}}
\newcommand{\ab}{\alpha ^{*}}

\begin{document}

\begin{center}
{\bf BIDIAGONAL PAIRS, THE LIE ALGEBRA $\SL$, AND THE QUANTUM GROUP $\uq$}
\end{center}

\medskip

\begin{center}
DARREN FUNK-NEUBAUER \\ {\it Department of Mathematics and Physics \\ Colorado State University - Pueblo \\ 2200 Bonforte Boulevard, Pueblo, CO 81001 USA \\ darren.funkneubauer@colostate-pueblo.edu}
\end{center}

\begin{abstract}
\noindent
We introduce a linear algebraic object called a bidiagonal pair.  Roughly speaking, a bidiagonal pair is a pair of diagonalizable linear transformations on a finite-dimensional vector space, each of which acts in a bidiagonal fashion on the eigenspaces of the other.  We associate to each bidiagonal pair a sequence of scalars called a parameter array.  Using this concept of a parameter array we present a classification of bidiagonal pairs up to isomorphism.  The statement of this classification does not explicitly mention the Lie algebra $\SL$ or the quantum group $\uq$.  However, its proof makes use of the finite-dimensional representation theory of $\SL$ and $\uq$.  In addition to the classification we make explicit the relationship between bidiagonal pairs and modules for $\SL$ and $\uq$.   \\

{\bf AMS classification code:}  Primary: 17B37; Secondary: 15A21, 15A30, 17B10, 81R50.

\end{abstract}

\section{Introduction}

\noindent
In this paper we introduce a type of linear algebraic object called a bidiagonal pair, and classify these objects.  Roughly speaking, a bidiagonal pair is a pair of diagonalizable linear transformations on a finite-dimensional vector space, each of which acts in a bidiagonal fashion on the eigenspaces of the other (see Definition \ref{def:bidiag} for the precise definition).  Associated to each bidiagonal pair is a sequence of scalars called a parameter array (see Definition \ref{def:parray}).  The main result of this paper is a classification of the bidiagonal pairs, up to isomorphism, using the concept of a parameter array (see Theorem \ref{thm:class}).  The statement of this classification does not make it clear how to {\it construct} the bidiagonal pairs in each isomorphism class.  Hence, we also present a construction of the bidiagonal pairs using the finite-dimensional modules for the Lie algebra $\SL$ and the quantum group $\uq$ (see Theorem \ref{thm:uqsl2bidiag} and Theorem \ref{thm:main3}).  The finite-dimensional modules for $\SL$ and $\uq$ are well known (see Lemmas \ref{thm:compred}, \ref{thm:sl2mods}, \ref{thm:uqcompred}, \ref{thm:uq2mods}).  \\  

\noindent
Bidiagonal pairs had their genesis in the study of the well-known quantum algebras $U_q(\widehat{\mathfrak{sl}}_2)$ and $\uq$.  Bidiagonal pairs arose from the discovery of some new presentations of these algebras, which are now referred to as the equitable presentations.  The equitable presentations came about through an attempt to classify a type of linear algebraic object called a tridiagonal pair.  See below for more information on the equitable presentations and tridiagonal pairs.  Thus, the importance of bidiagonal pairs lies in the fact that they provide insight into the relationships between several closely connected algebraic objects.  Although this paper is the first to explicitly define bidiagonal pairs, they appear implicitly in \cite{Benkart04, Funk-Neubauer07, Funk-Neubauer09, ItoTer073, ItoTerWang06}.  Bidiagonal pairs were used in \cite{ItoTer073} to construct irreducible $U_q(\widehat{\mathfrak{sl}}_2)$-modules starting from a certain type of tridiagonal pair.  In \cite{Funk-Neubauer09} bidiagonal pairs were involved in using a certain type of tridiagonal pair to construct irreducible modules for the $q$-tetrahedron algebra.  See below for more information on the $q$-tetrahedron algebra.  Bidiagonal pairs were used in \cite{Benkart04} to construct irreducible $U_q(\widehat{\mathfrak{sl}}_2)$-modules starting from irreducible modules for the Borel subalgebra of $U_q(\widehat{\mathfrak{sl}}_2)$.  In \cite{Funk-Neubauer07} bidiagonal pairs were a part of the construction of $U_q(\widehat{\mathfrak{sl}}_2)$-modules from certain raising and lowering maps satisfying the $q$-Serre relations.  Bidiagonal pairs were used in \cite{ItoTerWang06} to show that the generators from the equitable presentation of $\uq$ have invertible actions on each finite-dimensional $\uq$-module.         \\

\noindent
We now discuss the equitable basis for $\SL$ and the equitable presentation of $\uq$ because they will be used in our construction of bidiagonal pairs. \\
 
\noindent
The equitable basis for $\SL$ was first introduced in \cite{HarTer07} as part of the study of the tetrahedron Lie algebra \cite[Definition 1.1]{HarTer07}.  The tetrahedron  algebra has been used in the ongoing investigation of tridiagonal pairs, and is closely related to a number of well-known Lie algebras including the Onsager algebra \cite[Proposition 4.7]{HarTer07}, the $\SL$ loop algebra \cite[Proposition 5.7]{HarTer07}, and the three point $\SL$ loop algebra \cite[Proposition 6.5]{HarTer07}.  For more information on the tetrahedron algebra, see \cite{Benkart07, Elduque07, Hartwig07, ItoTerinpress4}.  For more information on the equitable basis for $\SL$, see \cite{Al-Najjar10, Benkart10}.  \\

\noindent
The equitable presentation of $\uq$ was first introduced in \cite{ItoTerWang06}.  This presentation of $\uq$ was used in the study of the $q$-tetrahedron algebra \cite{ItoTer072}.  The $q$-tetrahedron algebra has been used in the ongoing investigation of tridiagonal pairs, and is closely related to a number of well-known algebras including the $U_q(\mathfrak{sl}_2)$ loop algebra \cite[Proposition 8.3]{ItoTer072}, and positive part of $U_q(\widehat{\mathfrak{sl}}_2)$ \cite[Proposition 9.4]{ItoTer072}.  For more information on the $q$-tetrahedron algebra, see \cite{Funk-Neubauer09, ItoTer071, ItoTerinpress1, ItoTerinpress2, Miki10}.  We note that there exists an equitable presentation for the quantum group $U_q(\mathfrak{g})$, where $\mathfrak{g}$ is any symmetrizable Kac-Moody algebra \cite{Terinpress}. \\   

\noindent
We now offer some additional information on tridiagonal pairs because of their close connection to bidiagonal pairs.  The precise definition of a tridiagonal pair is given in \cite[Definition 1.1]{ItoTanTer01}.  Tridiagonal pairs originally arose in algebraic combinatorics through the study of a combinatorial object called a P- and Q-polynomial association scheme \cite{Bannai84, ItoTanTer01, Leonard82, Ter92}.  Since then they have appeared in many other areas of mathematics.  For instance, tridiagonal pairs arise in representation theory \cite{Al-Najjar05, Benkart04, Funk-Neubauer09, Hartwig07, ItoTer071, ItoTer072, ItoTer073, ItoTer074, ItoTer075, ItoTerinpress4, ItoTerinpress0, Koelink96, Koelink00, Koornwinder93, Noumi92, Rosengren99, Rosengren07, Ter012, TerVid04}, orthogonal polynomials and special functions \cite{Ter011, Ter041, Ter042}, partially ordered sets \cite{Ter90, Ter03}, statistical mechanics \cite{Baseilhac05, Baseilhac052, Baseilhac06, HarTer07}, and classical mechanics \cite{Zhedanov02}.  A certain special case of a tridiagonal pair, called a Leonard pair, has also been the subject of much research.  For example, the Leonard pairs were classified in \cite{Ter011, Ter051}.  They correspond to a family of orthogonal polynomials consisting of the $q$-Racah polynomials, and related polynomials in the terminating branch of the Askey scheme \cite{Askey79, Koekoek98, Ter041}.  For further information on tridiagonal pairs and Leonard pairs, see \cite{Al-Najjar11, Al-Najjar04, Al-Najjarinpress, Curtin071, Curtin072, Hartwig05, ItoTer04, ItoTer09, ItoTer105, Nomura051, Nomura052, Nomura053, Nomura061,  Nomura062, Nomura063, Nomura071, Nomura072, Nomura073, Nomura074, Nomurainpress1, Nomurainpress2, Nomura081, Nomurainpress3, Nomura09, Nomura092, Nomura102, Nomura10, Nomura11, Ter02, Ter052, Ter053, Vidarinpress, Vid07, Vid08}. \\

\section{Bidiagonal pairs}

\noindent
In this section we present the definition of a bidiagonal pair, and make some observations about this definition. \\

\noindent
We present the following lemma in order to motivate the definition of a bidiagonal pair.

\begin{lemma}
\label{thm:raise}
Let $\K$ denote a field.  Let $V$ denote a vector space over $\K$ with finite positive dimension.  Let $A:V\to V$ and $\Ab :V \to V$ denote linear transformations, and define $[A, A^*] := AA^* - A^*A$.  Suppose that there exists an ordering $\lbrace V_i \rbrace _{i=0}^d$ (resp.~$\lbrace \Vb_i \rbrace _{i=0}^{\delta}$) of the eigenspaces of $A$ (resp.~$\Ab$) with 
\begin{eqnarray*}
A \Vb_i \subseteq \Vb_i + \Vb_{i+1} \qquad (0 \leq i \leq \delta), \\
\Ab V_i \subseteq V_i + V_{i+1} \qquad (0 \leq i \leq d), 
\end{eqnarray*}
where $\Vb_{\delta+1} = 0$ and $V_{d+1} = 0$.  Then the following {\rm(i)},{\rm(ii)} hold.
\begin{enumerate}
\item[\rm (i)] For $0 \leq i \leq \delta$, the restriction $[A, A^*] |_{V_i^*}$  maps $V^*_i$ into $V_{i+1}^*$. 
\item[\rm (ii)] For $0 \leq i \leq d$, the restriction $[A, A^*] |_{V_i}$  maps $V_i$ into $V_{i+1}$.
\end{enumerate}
\end{lemma}

\noindent
{\it Proof:}  (i) Let $\theta_i^*$ (resp.~$\theta_{i+1}^*$) denote the eigenvalue of $A^*$ corresponding to $V_i^*$ (resp.~$V_{i+1}^*$).  Let $x \in V_i^*$.  By our assumption, there exists $y \in V_i^*$ and $z \in V_{i+1}^*$ such that $Ax = y + z$.  Then
\begin{eqnarray*}
[A, \Ab]x
&=& ( \theta_i^* - A^*) Ax \\
&=& ( \theta_i^* - A^*) y + ( \theta_i^* - A^*) z \\
&=& ( \theta_i^* - \theta_{i+1}^*) z.
\end{eqnarray*}
Therefore, $[A, \Ab]x \in V_{i+1}^*$, and so the restriction $[A, A^*] |_{V_i^*}$  maps $V^*_i$ into $V_{i+1}^*$.  \\
(ii) Similar to (i).  
\hfill $\Box$ \\

\begin{definition}
\rm
\label{def:bidiag}
Let $\K$ denote a field.  Let $V$ denote a vector space over $\K$ with finite positive dimension.  By a {\it bidiagonal pair on $V$},
we mean an ordered pair of linear transformations $A:V\to V$ and $\Ab :V \to V$ that satisfy the following three conditions.
\begin{enumerate}
\item[\rm (i)] Each of $A,\, \Ab$ is diagonalizable.
\item[\rm (ii)] There exists an ordering $\lbrace V_i \rbrace _{i=0}^d$ (resp.~$\lbrace \Vb_i \rbrace _{i=0}^{\delta}$) of the eigenspaces of $A$ (resp.~$\Ab$) with 
\begin{align}
\label{eq:b1}
A \Vb_i \subseteq \Vb_i + \Vb_{i+1} \qquad (0 \leq i \leq \delta), \\
\label{eq:b2}
\Ab V_i \subseteq V_i + V_{i+1} \qquad (0 \leq i \leq d),
\end{align}
where $\Vb_{\delta+1} = 0$ and $V_{d+1} = 0$.
\item[\rm (iii)] The restrictions
\begin{align}
\label{eq:b3}
[A, A^*]^{\delta-2i} |_{\Vb_i} : \Vb_i \rightarrow \Vb_{\delta-i} \qquad (0 \leq i \leq \delta/2), \\
\label{eq:b4}
[A, A^*]^{d-2i} |_{V_i} : V_i \rightarrow V_{d-i} \qquad (0 \leq i \leq d/2), 
\end{align}
are bijections, where $[A, A^*] := AA^* - A^*A$.
\end{enumerate}
\end{definition}

\noindent
We call $V$ the {\it vector space underlying $A, \, A^*$}, and say $A, \, A^*$ is {\it over} $\K$.  We call $d$ the {\it diameter} of $A, \, A^*$.

\begin{note}
\label{*notation}
\rm
According to a common notational convention, $A^*$ denotes the conjugate transpose of $A$.  We are {\it not} using this convention. In a bidiagonal pair the linear transformations $A$ and $A^*$ are arbitrary subject to (i)--(iii) above.
\end{note}

\begin{lemma}
\label{thm:dequalsdelta}
With reference to the notation in Definition \ref{def:bidiag}, we have $d=\delta$.
\end{lemma}

\noindent
The proof of Lemma \ref{thm:dequalsdelta} is contained in the appendix of this paper.  In view of Lemma \ref{thm:dequalsdelta}, for the remainder of the paper we use $d$ to index the eigenspaces of both $A$ and $A^*$.  \\

\noindent
For the remainder of this paper, we assume that the field $\K$ in Definition \ref{def:bidiag} is algebraically closed and characteristic zero.  

\begin{lemma}
\label{thm:injsurj}
With reference to Definition \ref{def:bidiag}, for $0 \leq i < d/2$ (resp.~$d/2 \leq i \leq d$), and for $0 \leq j \leq d-2i$ (resp.~$0 \leq j \leq d-i$), the restrictions
\begin{align*}
[A, \Ab]^j |_{\Vb_i} : \Vb_i \rightarrow \Vb_{i+j}, \\
[A, \Ab]^j |_{V_i} : V_i \rightarrow V_{i+j}
\end{align*}
are injections (resp.~surjections).
\end{lemma}

\noindent
{\it Proof:} Immediate from Definition \ref{def:bidiag}(iii).
\hfill $\Box $ \\

\noindent
An ordering of the eigenspaces of $A$ (resp.~$A^*$) is called {\it standard} whenever this ordering satisfies (\ref{eq:b2}) (resp.~(\ref{eq:b1})).  If  $\lbrace V_i \rbrace _{i=0}^d$ (resp.~$\lbrace \Vb_i \rbrace _{i=0}^d$) is a standard ordering of the eigenspaces of $A$ (resp.~$\Ab$) then no other ordering of $\lbrace V_i \rbrace _{i=0}^d$ (resp.~$\lbrace \Vb_i \rbrace _{i=0}^d$) is standard.  Thus, the bidiagonal pair $A, \, A^*$ uniquely determines the standard orderings of the eigenspaces of $A$ and $A^*$.  Let $\lbrace V_i \rbrace _{i=0}^d$ (resp.~$\lbrace \Vb_i \rbrace _{i=0}^d$) denote the standard ordering of the eigenspaces of $A$ (resp.~$\Ab$).  For $0 \leq i \leq d$, let $\theta_i$ (resp.~$\Tb_i $) denote the eigenvalue of $A$ (resp.~$\Ab$) corresponding to $V_i$ (resp.~$\Vb_i$).  We call $\lbrace \theta_i \rbrace _{i=0}^d$ (resp.~$\lbrace \Tb_i \rbrace _{i=0}^d$) the {\it eigenvalue} (resp.~{\it dual eigenvalue}) {\it sequence} of $A, \, \Ab$. \\

\noindent
For the remainder of this section, we adopt the following assumption.  Let $A, \, \Ab$ denote a bidiagonal pair on $V$, and let $\lbrace V_i \rbrace _{i=0}^d$ (resp.~$\lbrace V^*_i \rbrace _{i=0}^d$) denote the standard ordering of the eigenspaces of $A$ (resp.~$\Ab$).  Let $\lbrace \theta_i \rbrace _{i=0}^d$ (resp.~$\lbrace \Tb_i \rbrace _{i=0}^d$) denote the eigenvalue (resp.~dual eigenvalue) sequence of $A, \, \Ab$.

\begin{lemma}
\label{thm:olddefn}
For $0 \leq i \leq d$, we have
\begin{align}
\label{olddefn1}
(A - \theta_i I) \Vb_i \subseteq \Vb_{i+1}, \\
\label{olddefn2}
(\Ab - \theta^*_i I) V_i \subseteq V_{i+1}.
\end{align}
\end{lemma}

\noindent
We postpone the proof of Lemma \ref{thm:olddefn} until Section 6.

\begin{lemma}
\label{thm:samedim}
For $0 \leq i \leq d$, the spaces $V_i$, $V_{d-i}$, $\Vb_i$, $\Vb_{d-i}$ all have the same dimension.
\end{lemma}

\noindent
{\it Proof:}  Using (\ref{olddefn2}) we see that, with respect to an appropriate basis for $V$, $A^*$ can be represented as a lower triangular matrix with diagonal entries $\Tb_0, \, \Tb_1, \, \ldots, \, \Tb_d$, where $\Tb_i$ appears $\dim(V_i)$ times.  Hence, the multiplicity of $\Tb_i$ as a root of the characteristic polynomial of $\Ab$ is $\dim(V_i)$.  However, since $A^*$ is diagonalizable, we find that the multiplicity of $\Tb_i$ as a root of the characteristic polynomial of $\Ab$ is also $\dim(V^*_i)$.  Therefore, $\dim(V_i) = \dim(V^*_i)$.  Combining this with Definition \ref{def:bidiag}(iii) we obtain the result.
\hfill $\Box$ \\

\begin{definition}
\rm
\label{def:shape}
With reference to Lemma \ref{thm:samedim}, for $0 \leq i \leq d$, let $\rho_i$ denote the common dimension of $V_i$, $V_{d-i}$, $\Vb_i$, $\Vb_{d-i}$.  We refer to the sequence $\lbrace \rho_i \rbrace _{i=0}^d$ as the {\it shape} of $A, \, A^*$.
\end{definition}

\noindent
We now define the notion of the parameter array of $A, \, A^*$ which will be used in stating our main result (see Theorem \ref{thm:class}).

\begin{definition}
\rm
\label{def:parray}
By the {\it parameter array} of $A, \, A^*$, we mean the sequence \\ $(\lbrace \theta_i \rbrace _{i=0}^d; \lbrace \theta^*_i \rbrace _{i=0}^d; \lbrace \rho_i \rbrace _{i=0}^d)$, where $\lbrace \theta_i \rbrace _{i=0}^d$ (resp.~$\lbrace \theta^*_i \rbrace _{i=0}^d$) is the eigenvalue (resp.~dual eigenvalue) sequence of $A, \, A^*$, and $\lbrace \rho_i \rbrace _{i=0}^d$ is the shape of $A, \, A^*$.
\end{definition}

\begin{lemma}
\label{affineisbidiag}
Let $(\lbrace \theta_i \rbrace _{i=0}^d; \lbrace \theta^*_i \rbrace _{i=0}^d; \lbrace \rho_i \rbrace _{i=0}^d)$ denote the parameter array of $A,\, A^*$.  Let $p, \, q, \, r, \, s$ denote scalars in $\K$ with both $p, \, r$ nonzero.   Then $pA + qI , \, rA^* + sI$ is a bidiagonal pair on $V$ with parameter array $(\lbrace p \theta_i + q \rbrace _{i=0}^d; \lbrace r \theta^*_i  + s \rbrace _{i=0}^d; \lbrace \rho_i \rbrace _{i=0}^d)$.
\end{lemma}

\noindent
{\it Proof:}  For $v \in V$ and $0 \leq i \leq d$, we have $Av = \theta_i v$ if and only if $(pA + qI)v = (p \theta_i + q)v$.  Thus, $p \theta_i + q$ is an eigenvalue of $pA + qI$, and $V_i$ is the eigenspace of $pA + qI$ corresponding to $p \theta_i + q$.  We know $V = \sum_{i=0}^d V_i$ (direct sum) since $A$ is diagonalizable.  So $pA + qI$ is diagonalizable.  Similarly, $r \Tb_i + s$ is an eigenvalue of $rA^* + sI$, $\Vb_i$ is the eigenspace of $rA^* + sI$ corresponding to $r \Tb_i + s$, and $rA^* + sI$ is diagonalizable.  It is immediate from (\ref{eq:b1}) (resp.~(\ref{eq:b2})) that $(pA + qI) \Vb_i \subseteq \Vb_i + \Vb_{i+1}$ (resp.~$(rA^* + sI) V_i \subseteq V_i + V_{i+1}$).  Observe $[pA + qI, rA^* + sI] = pr [A, A^*]$.  From this, and (\ref{eq:b3}) (resp.~(\ref{eq:b4})) we have $[pA + qI, rA^* + sI]^{d-2i} |_{\Vb_i} : \Vb_i \rightarrow \Vb_{d-i}$ (resp.~$[pA + qI, rA^* + sI]^{d-2i} |_{V_i} : V_i \rightarrow V_{d-i}$) is a  bijection.  We have now shown that $pA + qI , rA^* + sI$ is a bidiagonal pair on $V$, and that $\lbrace V_i \rbrace _{i=0}^d$ (resp.~$\lbrace V^*_i \rbrace _{i=0}^d$) is a standard ordering of the eigenspaces of $pA + qI$ (resp.~$rA^* + sI$).  Thus, $pA + qI , rA^* + sI$ has eigenvalue sequence $\lbrace p \theta_i + q \rbrace _{i=0}^d$, dual eigenvalue sequence $\lbrace r \theta^*_i  + s \rbrace _{i=0}^d$, and shape $\lbrace \rho_i \rbrace _{i=0}^d$.
\hfill $\Box$ \\

\noindent
We now define the notion of isomorphism for bidiagonal pairs.

\begin{definition}
\rm
\label{def:isom}
Let $A, \, A ^*$ and $B, \, B^*$ denote bidiagonal pairs over $\K$.  Let $V$ (resp.~$\widetilde{V}$) denote the vector space underlying $A, \, A ^*$ (resp.~$B, \, B^*$).  By an {\it isomorphism of bidiagonal pairs from $A, \, A ^*$ to $B, \, B^*$}, we mean a vector space isomorphism $\mu : V \rightarrow \widetilde{V}$ such that $\mu A = B \mu$ and $\mu A ^* = B ^* \mu$.  We say $A, \, A ^*$ and $B, B ^*$ are {\it isomorphic} whenever there exists an isomorphism of bidiagonal pairs from $A, \, A ^*$ to $B, B ^*$.
\end{definition}

\begin{lemma}
\label{thm:isompa}
Let $A, \, A ^*$ and $B, \, B^*$ denote bidiagonal pairs over $\K$.  Then  $A, \, A ^*$ and $B, B ^*$ are {\it isomorphic} if and only if the parameter array of $A, \, A^*$ equals the parameter array of $B, \, B^*$.
\end{lemma}

\noindent
We postpone the proof of Lemma \ref{thm:isompa} until Section 13. \\

\noindent
We now define the notion of affine equivalence which will be used in stating one of our main theorems (see Theorem \ref{thm:main2}). 

\begin{definition}
\rm
\label{def:affine}
Let $A, \, A ^*$ and $B, \, B^*$ denote bidiagonal pairs over $\K$.  We say $A, \, A ^*$ is {\it affine equivalent} to $B, B ^*$ if there exist scalars $p,\, q,\, r, \, s \in \K$ with both $p, \, r$ nonzero such that $A, \, \Ab$ and $p B + qI, \, r B^* + sI$ are isomorphic.
\end{definition}

\noindent
The notion of affine equivalence is an equivalence relation on bidiagonal pairs.  If two bidiagonal pairs are isomorphic then they are in the same affine equivalence class.

\section{The Lie algebra $\SL$}

\noindent
In this section we recall $\SL$ and its finite-dimensional modules.

\begin{definition}
\rm
\label{def:usl2}
Let $\SL$ denote the Lie algebra over $\K$ that has a basis $h, \, e, \, f$ and Lie bracket
\begin{eqnarray*}
[h,e] = 2e, \qquad
[h,f] = -2f, \qquad
[e,f] = h.
\end{eqnarray*}
\end{definition}

\begin{theorem}
\cite[Lemma 3.2]{HarTer07}
\label{thm:ueq}
The Lie algebra $\SL$ is isomorphic to the Lie algebra over $\K$ that has basis $X,\, Y, \, Z$ and Lie bracket
\begin{eqnarray}
\label{equitablepres}
[X,Y] = 2X + 2Y, \qquad
[Y,Z] = 2Y + 2Z, \qquad
[Z,X] =2Z + 2X.
\end{eqnarray}
An isomorphism with the presentation in Definition \ref{def:usl2} is given by:
\begin{eqnarray*}
X &\rightarrow& 2e - h, \\
Y &\rightarrow& -2f - h, \\
Z &\rightarrow& h.
\end{eqnarray*}
The inverse of this isomorphism is given by:
\begin{eqnarray*}
e &\rightarrow& (X + Z)/2, \\
f &\rightarrow& -(Y + Z)/2, \\
h &\rightarrow& Z.
\end{eqnarray*}
\end{theorem}

\noindent
By an {\it equitable basis} for $\SL$, we mean a basis $X, \, Y, \, Z$ that satisfies (\ref{equitablepres}).

\begin{definition}
\rm
Given an ordered pair $Y, \, Z$ of elements in $\SL$, we call this pair {\it equitable} whenever there exists an element $X$ in $\SL$ such that $X, \, Y, \, Z$ is an equitable basis for $\SL$.
\label{def:equitpair}

\end{definition}

\begin{note}
\label{def:psi}
\rm
From Theorem \ref{thm:ueq} it is clear that there is a Lie algebra automorphism $\psi$ of $\SL$ with order $3$ such that
\begin{eqnarray*}
\psi(X) = Y, \qquad \psi(Y) = Z, \qquad \psi(Z) = X.
\end{eqnarray*}
\end{note}

\noindent
The following two lemmas give a description of all finite-dimensional $\SL$-modules.

\begin{lemma}
\cite[Theorem 6.3]{Humph72}
\label{thm:compred}
Each finite-dimensional $\SL$-module $V$ is completely reducible; this means that $V$ is a direct sum of irreducible $\SL$-modules.
\end{lemma}

\noindent
The finite-dimensional irreducible $\SL$-modules are described as follows.

\begin{lemma}
\cite[Theorem 7.2]{Humph72}
\label{thm:sl2mods}
There exists a family of finite-dimensional irreducible $\SL$-modules
\begin{eqnarray*}
V(d), \qquad  d=0,1,2,\ldots
\end{eqnarray*}
with the following properties:  $V(d)$ has a basis $\lbrace v_i \rbrace _{i=0}^d$ such that $h.v_i = (d-2i) v_i$ for $0 \leq i \leq d$, $f.v_i = (i+1) v_{i+1}$ for $0 \leq i \leq d$, where $v_{d+1}=0$, and $e.v_i = (d-i+1) v_{i-1}$ for $0 \leq i \leq d$, where $v_{-1}=0$.  Moreover, every finite-dimensional irreducible $\SL$-module is isomorphic to exactly one of the modules $V(d)$.
\end{lemma}

\begin{lemma}
\label{thm:2pieces}
Let $V$ denote an $\SL$-module with finite positive dimension (not necessarily irreducible).  Define 
\begin{eqnarray*}
V_{\hbox{even}} &:=& \operatorname{span}\{\,v \in V\, |\, h.v = i \, v, \,\, i \in \Z, \, i \, even \,\}, \\
V_{\hbox{odd}} &:=& \operatorname{span}\{\,v \in V\, |\, h.v = i \, v, \,\, i \in \Z, \, i \, odd \,\}.
\end{eqnarray*}
Then $V_{\hbox{even}}$ and $V_{\hbox{odd}}$ are $\SL$-modules, and $V = V_{\hbox{even}} + V_{\hbox{odd}}$ (direct sum).
\end{lemma}

\noindent
{\it Proof:}  By construction each of $V_{\hbox{even}}$ and $V_{\hbox{odd}}$ are invariant under the action of $h$.  Using $[h,e] = 2e$ and $[h,f] = -2f$, we find that each of $V_{\hbox{even}}$ and $V_{\hbox{odd}}$ is invariant under the action of $e$ and $f$.  Thus, the subspaces $V_{\hbox{even}}$ and $V_{\hbox{odd}}$ are both $\SL$-submodules of $V$.  Combining Lemma \ref{thm:compred} and Lemma \ref{thm:sl2mods}, we find that the action of $h$ on $V$ is diagonalizable, and all the eigenvalues of this action are integers.  So $V$ is the direct sum of eigenspaces for the action of $h$ on $V$, and the result follows.
\hfill $\Box $ \\

\noindent
The following definition will be used in stating two of our main theorems (see Theorem \ref{thm:uqsl2bidiag} and Theorem \ref{thm:main3}).

\begin{definition}
\rm
\label{def:pure}
Let $V$ denote an $\SL$-module with finite positive dimension.  With reference to Lemma \ref{thm:2pieces}, we say $V$ is {\it segregated} whenever $V=V_{\hbox{even}}$ or $V=V_{\hbox{odd}}$.
\end{definition}

\section{The quantum group $\uq$}

\noindent
In this section we recall $\uq$ and its finite-dimensional modules.  \\

\noindent
In this section we assume that $q$ is a nonzero scalar in $\K$ which is not a root of unity.  For each nonnegative integer $n$, define
\begin{eqnarray*}
[n] = \frac{q^{n} - q^{-n}}{q-q^{-1}}.
\end{eqnarray*}

\begin{definition}
\rm
\label{def:uqsl2}
Let $\uq$ denote the unital associative $\K$-algebra with generators \\ $k, \, k^{-1}, \, e, \, f$ and the following relations:
\begin{eqnarray*}
kk^{-1} &=& k^{-1}k = 1, \\
ke &=& q^2 ek, \\
kf &=& q^{-2}fk, \\
ef - fe &=& \frac{k - k^{-1}}{q - q^{-1}}.
\end{eqnarray*}
\end{definition}

\begin{theorem}
\cite[Theorem 2.1]{ItoTerWang06}
\label{thm:uqeq}
The algebra $\uq$ is isomorphic to the unital associative $\K$-algebra with generators $x, \, x^{-1}, \, y, \, z$ and the following relations:
\begin{eqnarray*}
xx^{-1} = x^{-1}x &=& 1, \\
\frac{qxy - q^{-1}yx}{q - q^{-1}} &=& 1, \\
\frac{qyz - q^{-1}zy}{q - q^{-1}} &=& 1, \\
\frac{qzx - q^{-1}xz}{q - q^{-1}} &=& 1.
\end{eqnarray*}
An isomorphism with the presentation in Definition \ref{def:uqsl2} is given by:
\begin{eqnarray*}
x^{\pm1} &\rightarrow& k^{\pm1}, \\
y &\rightarrow& k^{-1} + f(q-q^{-1}), \\
z &\rightarrow& k^{-1} - k^{-1}eq(q-q^{-1}).
\end{eqnarray*}
The inverse of this isomorphism is given by:
\begin{eqnarray*}
k^{\pm1} &\rightarrow& x^{\pm1}, \\
f &\rightarrow& (y - x^{-1})(q-q^{-1}), \\
e &\rightarrow& (1 - xz)q^{-1}(q-q^{-1}).
\end{eqnarray*}
\end{theorem}

\noindent
By the {\it equitable presentation} for $\uq$, we mean the presentation given in Theorem \ref{thm:uqeq}.  We call $x, \, x^{-1}, \, y, \, z$ the {\it equitable generators} for $\uq$.

\begin{definition}
\rm
\label{def:equitpair2}
Given an ordered pair $y, \, z$ of elements in $\uq$, we call this pair {\it equitable} whenever there exists elements $x, \, x^{-1}$ in $\uq$ such that $x, \, x^{-1}, \, y, \, z$ are equitable generators for $\uq$.
\end{definition}

\begin{lemma}
\cite[Theorem 7.5]{ItoTerWang06}
\label{thm:omega}
Let $x,\, y, \, z$ be the equitable generators for $\uq$.  There exists an invertible linear operator $\Omega$ that acts on finite-dimensional $\uq$-modules, and satisfies
\begin{eqnarray*}
\Omega ^{-1} \, x \, \Omega = y, \qquad \Omega ^{-1} \, y \, \Omega = z, \qquad \Omega ^{-1} \, z \, \Omega = x. \end{eqnarray*}
\end{lemma}

\noindent
The following two lemmas give a description of all finite-dimensional $\uq$-modules.

\begin{lemma}
\cite[Theorem 2.3, 2.9]{Jantzen96}
\label{thm:uqcompred}
Each finite-dimensional $\uq$-module $V$ is completely reducible; this means that $V$ is a direct sum of irreducible $\uq$-modules.
\end{lemma}

\noindent
The finite-dimensional irreducible $\uq$-modules are described as follows.

\begin{lemma}
\cite[Theorem 2.6]{Jantzen96}
\label{thm:uq2mods}
There exists a family of finite-dimensional irreducible $\uq$-modules
\begin{eqnarray*}
V(d, \epsilon), \qquad \epsilon \in \{1, -1 \}, \qquad  d=0,1,2,\ldots
\end{eqnarray*}
with the following properties:  $V(d, \epsilon)$ has a basis $\lbrace v_i \rbrace _{i=0}^d$ such that $k.v_i = \epsilon q^{d-2i} v_i$ for $0 \leq i \leq d$, $f.v_i = [i+1] v_{i+1}$ for $0 \leq i \leq d$, where $v_{d+1}=0$, and $e.v_i = \epsilon [d-i+1] v_{i-1}$ for $0 \leq i \leq d$, where $v_{-1}=0$.  Moreover, every finite-dimensional irreducible $\uq$-module is isomorphic to exactly one of the modules $V(d, \epsilon)$.
\end{lemma}

\begin{lemma}
\label{thm:4pieces}
Let $V$ denote a $\uq$-module with finite positive dimension (not necessarily irreducible).  For $\epsilon \in \{1, -1 \}$, define
\begin{eqnarray*}
V_{\hbox{even}}^{\epsilon} &:=& \operatorname{span}\{\,v \in V\, |\, k.v = \epsilon \, q^i \, v, \,\, i \in \Z, \, i \, even \,\}, \\
V_{\hbox{odd}}^{\epsilon} &:=& \operatorname{span}\{\,v \in V\, |\, k.v = \epsilon \, q^i \, v, \,\, i \in \Z, \, i \, odd \,\}.
\end{eqnarray*}
Then $V_{\hbox{even}}^{1}, \, V_{\hbox{even}}^{-1}, \, V_{\hbox{odd}}^{1}, \, V_{\hbox{odd}}^{-1}$ are $\uq$-modules, and $V = V_{\hbox{even}}^{1} + V_{\hbox{even}}^{-1} + V_{\hbox{odd}}^{1} + V_{\hbox{odd}}^{-1}$ (direct sum).
\end{lemma}

\noindent
{\it Proof:}  Let $\epsilon \in \{1, -1 \}$.  By construction each of $V_{\hbox{even}}^{\epsilon}$ and $V_{\hbox{odd}}^{\epsilon}$ are invariant under the action of $k$ and $k^{-1}$.  Using $ke = q^2 ek$ and $kf = q^{-2} fk$, we find that each of $V_{\hbox{even}}^{\epsilon}$ and $V_{\hbox{odd}}^{\epsilon}$ are invariant under the action of $e$ and $f$.  Thus, the subspaces $V_{\hbox{even}}^{\epsilon}$ and $V_{\hbox{odd}}^{\epsilon}$ are both $\uq$-submodules of $V$.  Combining Lemma \ref{thm:uqcompred} and Lemma \ref{thm:uq2mods}, we find that the action of $k$ on $V$ is diagonalizable, and all the eigenvalues of this action are plus or minus integer powers of $q$.  So $V$ is the direct sum of eigenspaces for the action of $k$ on $V$, and the result follows.
\hfill $\Box $ \\

\noindent
The following definition will be used in stating two of our main theorems (see Theorem \ref{thm:uqsl2bidiag} and Theorem \ref{thm:main3}).

\begin{definition}
\rm
\label{def:uqpure}
Let $V$ denote a $\uq$-module with finite positive dimension.  With reference to Lemma \ref{thm:4pieces}, we say  $V$ is {\it segregated} whenever $V=V_{\hbox{even}}^1$ or $V=V_{\hbox{odd}}^1$.
\end{definition}

\section{The main theorems}

\noindent
The five theorems in this section make up the main conclusions of the paper.  The following theorem provides a classification of bidiagonal pairs up to isomorphism.

\begin{theorem}
\label{thm:class}
Let $d$ denote a nonnegative integer, and let
\begin{eqnarray}
\label{parameterarray}
(\lbrace \theta_i \rbrace _{i=0}^d; \lbrace \theta^*_i \rbrace _{i=0}^d; \lbrace \rho_i \rbrace _{i=0}^d)
\end{eqnarray}
denote a sequence of scalars taken from $\K$.  Then there exists a bidiagonal pair $A, \, A^*$ over $\K$ with parameter array (\ref{parameterarray}) if and only if {\rm (i)--(v)} below hold.
\begin{enumerate}
\item[\rm (i)] $\theta_i \neq \theta_j$, $\theta^*_i \neq \theta^*_j$ for $0 \leq i,j \leq d$ and $i \neq j$.
\item[\rm (ii)] The expressions
\begin{eqnarray*}
\frac{\theta_{i+1} - \theta_i}{\theta_i - \theta_{i-1}}, \qquad \frac{\Tb_i - \Tb_{i-1}}{\Tb_{i+1} - \Tb_i}
\end{eqnarray*}
are equal and independent of $i$ for $1 \leq i \leq d-1$.
\item[\rm (iii)] $\rho_i$ is a positive integer for $0 \leq i \leq d$.
\item[\rm (iv)] $\rho_i = \rho_{d-i}$ for $0 \leq i \leq d$.
\item[\rm (v)] $\rho_{i} \leq \rho_{i+1}$ for $0 \leq i < d/2$.
\end{enumerate}
Suppose that {\rm(i)--(v)} hold.  Then $A, \, A^*$ is unique up to isomorphism of bidiagonal pairs.
\end{theorem}

\noindent
We refer to Theorem \ref{thm:class} as the classification theorem. 

\begin{note}
\label{vac}
\rm
With reference to Theorem \ref{thm:class}, for $d=0$ we regard conditions (i), (ii), (v) as holding since they are vacuously true.  Similarly, for $d=1$ we regard condition (ii) as holding since it is vacuously true.
\end{note}

\noindent
It is not clear from the statement of Theorem \ref{thm:class} how bidiagonal pairs are related to finite-dimensional modules for $\SL$ and $\uq$.  Also, it is not clear from Theorem \ref{thm:class} how the bidiagonal pairs in each isomorphism class are constructed.  The next four theorems address these issues.

\begin{theorem}
\label{thm:main1}
Let $A, \, \Ab$ denote a bidiagonal pair.  Then there exists a sequence of scalars $b, \, \alpha, \, \ab, \, \gamma$ in $\K$ with $b$ nonzero such that
\begin{eqnarray}
\label{main1.1}
A \Ab - b \Ab A - \alpha A - \ab \Ab - \gamma I = 0.
\end{eqnarray}
This sequence of scalars is uniquely determined by the pair $A, \, \Ab$ provided $d \geq 2$.
\end{theorem}

\noindent
We refer to (\ref{main1.1}) as the fundamental bidiagonal relation. 

\begin{note}
\label{scalars}
\rm
See Lemma \ref{thm:relrec} for a description of how the scalars $b, \, \alpha, \, \alpha ^*, \, \gamma$ are related to the eigenvalues of $A$ and $A^*$.
\end{note}

\begin{definition}
\label{def:base}
\rm
We call the scalar $b$ from (\ref{main1.1}) the {\it base} of the bidiagonal pair $A, \, A^*$.  By Theorem \ref{thm:main1}, the base of $A, \, \Ab$ is well defined for $d \geq 2$.  For $d \leq 1$, (\ref{main1.1}) can be shown to hold when $b$ is any nonzero scalar in $\K$ (see the proof of Theorem \ref{thm:main1} in Section 8).  We adopt the convention that $b=1$ for $d \leq 1$.  With this convention, the base is well defined for each nonnegative integer $d$.
\end{definition}

\begin{note}
\label{bnot1}
\rm
Let $b$ denote the base of $A, \, A^*$.  Whenever $b \neq 1$, let $q \in \K$ denote a root of the polynomial $\lambda ^2 - b^{-1} \in \K[\lambda]$.  The scalar $q$ exists since $\K$ is algebraically closed.  Observe $q \neq 0$, and so $b = q^{-2}$.  By construction $b$ uniquely determines $q$ up to sign.
\end{note}

\noindent
The following definition will be used to state the next three theorems.

\begin{definition}
\label{def:reduced}
\rm
Let $A, \, A^*$ denote a bidiagonal pair of diameter $d$, and let $b$ denote its base.  We say $A, \, \Ab$ is {\it reduced} if either (i) or (ii) holds.
\begin{enumerate}
\item[\rm (i)] $b=1$ and the eigenvalue (resp.~dual eigenvalue) sequence of $A, \, \Ab$ is \\ $\lbrace 2i-d \rbrace _{i=0}^d$ (resp.~$\lbrace d-2i \rbrace _{i=0}^d$).
\item[\rm (ii)] $b \neq 1$ and the eigenvalue (resp.~dual eigenvalue) sequence of $A, \, \Ab$ is $\lbrace q^{d-2i} \rbrace _{i=0}^d$ (resp.~$\lbrace q^{2i-d} \rbrace _{i=0}^d$). \\
\end{enumerate}
\end{definition}

\begin{theorem}
\label{thm:main2}
Every bidiagonal pair is affine equivalent to a reduced bidiagonal pair.
\end{theorem}

\noindent
We refer to the result in Theorem \ref{thm:main2} as the reducibility of bidiagonal pairs. 

\begin{note}
\label{reducedcases}
\rm
We make the following observation in order to motivate the next two theorems.  Let $A, \, A^*$ denote a reduced bidiagonal pair, and let $b$ denote its base.  It is shown in Section 12 that for $b=1$, (\ref{main1.1}) takes the form
\begin{eqnarray*}
AA^* - A^* A - 2A - 2 A^* = 0.
\end{eqnarray*}
It is also shown in Section 12 that for $b \neq 1$, (\ref{main1.1}) takes the form
\begin{eqnarray*}
qAA^* - q^{-1} A^* A - (q - q^{-1}) I =  0.
\end{eqnarray*}
\end{note}

\noindent
The following two theorems explain the connection between reduced bidiagonal pairs and finite-dimensional modules for $\SL$ and $\uq$.

\begin{theorem}
\label{thm:uqsl2bidiag}
Let $V$ denote a segregated $\SL$-module (resp.~segregated $\uq$-module).  Then each equitable pair in $\SL$ (resp.~$\uq$) acts on $V$ as a reduced bidiagonal pair with base $1$ (resp.~base not equal to $1$).
\end{theorem}

\noindent
We refer to the result in Theorem \ref{thm:uqsl2bidiag} as equitable pairs act as bidiagonal pairs. \\

\noindent
The next theorem can be thought of as the converse of Theorem \ref{thm:uqsl2bidiag}.

\begin{theorem}
\label{thm:main3}
Let $A, \, \Ab$ denote a reduced bidiagonal pair on $V$, and let $b$ denote its base.  Then the following {\rm(i)},{\rm(ii)} hold.
\begin{enumerate}
\item[\rm (i)] Suppose that $b = 1$.  Let $Y, \, Z$ denote an equitable pair in $\SL$.  Then there exists a $\SL$-module structure on $V$ such that $(Y - A)V=0$ and $(Z - \Ab)V=0$.  The action of $X$ on $V$ is uniquely determined by $A, \, A^*$. This $\SL$-module structure on $V$ is segregated.
\item[\rm (ii)] Suppose that $b \neq 1$.  Let $y, \, z$ denote an equitable pair in $\uq$.  Then there exists a $\uq$-module structure on $V$ such that $(y - A)V=0$ and $(z - \Ab)V=0$.  The action of $x^{\pm1}$ on $V$ is uniquely determined by $A, \, A^*$.  This $\uq$-module structure on $V$ is segregated.
\end{enumerate}
\end{theorem}

\noindent
We refer to the result in Theorem \ref{thm:main3} as bidiagonal pairs act as equitable pairs. \\

\noindent
Theorem \ref{thm:main2} and Theorem \ref{thm:main3} together show that every bidiagonal pair is affine equivalent to a bidiagonal pair of the type constructed in Theorem \ref{thm:uqsl2bidiag}.

\section{Preliminaries}

\noindent
In this section we develop some more properties of bidiagonal pairs.  We also prove Lemma \ref{thm:olddefn}, namely that the eigenspaces of one transformation in a bidiagonal pair are raised by the other transformation. \\

\noindent
Throughout the remainder of the paper we will refer to the following assumption.

\begin{assumption}
\label{assume}
\rm
Let $V$ denote a vector space over $\K$ with finite positive dimension, and let $A,\Ab$ denote a bidiagonal pair on $V$.  Let $\lbrace V_i \rbrace _{i=0}^d$ (resp.~$\lbrace V^*_i \rbrace _{i=0}^d$) denote a standard ordering of the eigenspaces of $A$ (resp.~$\Ab$).  Let $(\lbrace \theta_i \rbrace _{i=0}^d; \lbrace \theta^*_i \rbrace _{i=0}^d; \lbrace \rho_i \rbrace _{i=0}^d)$ denote the parameter array of $A, \, \Ab$.
\end{assumption}

\begin{note}
\label{A^*A}
\rm
With reference to Assumption \ref{assume}, it is immediate from Definition \ref{def:bidiag} that $A^*, \, A$ is a bidiagonal pair on V with parameter array  $(\lbrace \theta_i^* \rbrace _{i=0}^d; \lbrace \theta_i \rbrace _{i=0}^d; \lbrace \rho_i \rbrace _{i=0}^d)$.
\end{note}

\begin{definition}
\label{def:decomp}
\rm
Let $V$ denote a vector space over $\K$ with finite positive dimension.  By a \\{\it decomposition} of $V$, we mean a sequence  $\lbrace U_i \rbrace _{i=0}^d$ consisting of nonzero subspaces of $V$ such that $V=\sum_{i=0}^{d} U_i$ (direct sum).  For notational convenience, let $U_{-1}:=0$ and $U_{d+1}:=0$.
\end{definition}

\noindent
Referring to Assumption \ref{assume}, the sequences $\lbrace V_i \rbrace _{i=0}^d$, $\lbrace V^*_i \rbrace _{i=0}^d$ are both decompositions of $V$. \\

\noindent
The following definition will be used throughout the paper.

\begin{definition}
\label{def:H_i}
\rm
With reference to Assumption \ref{assume}, define
\begin{align}
\label{H_i1}
H_i &:= \{\,v \in V_i\, |\, [A, A^*]^{d-2i+1}v=0\,\} \qquad (0 \leq i \leq d/2), \\
\label{H_i2}
H^*_i &:= \{\,v \in \Vb_i\, |\, [A, A^*]^{d-2i+1}v=0\,\} \qquad (0 \leq i \leq d/2).
\end{align}
\end{definition}

\noindent
The following two lemmas will be used in proving Lemma \ref{thm:olddefn}.

\begin{lemma}
\label{thm:refineVbi}
With reference to Assumption \ref{assume} and Definition \ref{def:H_i}, the following hold.
\begin{eqnarray}
\label{refineV1}
V_i &=\sum_{j=0}^{\min(i,d-i)} [A,A^*]^{i-j} H_j \qquad \hbox{(direct sum)} \qquad (0 \leq i \leq d),  \\
\label{refineV2}
\Vb_i &=\sum_{j=0}^{\min(i,d-i)} [A,A^*]^{i-j} H^*_j \qquad \hbox{(direct sum)} \qquad (0 \leq i \leq d).
\end{eqnarray}
\end{lemma}

\noindent
{\it Proof:}  First we prove (\ref{refineV1}).  \\
Case 1: $0 \leq i \leq d/2$.  The proof is by induction on $i$.  Observe the result holds for $i=0$, since $V_0=H_0$.  Next assume that $i \geq 1$. By induction and Lemma \ref{thm:injsurj}, we find
\begin{eqnarray}
\label{refineVbi1}
[A, A^*] V_{i-1}=\sum_{j=0}^{i-1} [A, A^*]^{i-j} H_j \qquad  \hbox{(direct sum)}.
\end{eqnarray}
We now show
\begin{eqnarray}
\label{refineVbi2}
V_i = [A, A^*] V_{i-1} + H_i \qquad \hbox{(direct sum)}.
\end{eqnarray}
Lemma \ref{thm:raise} and (\ref{H_i1}) yield $[A, A^*] V_{i-1} + H_i \subseteq V_i$.  We now show that $V_i \subseteq [A, A^*] V_{i-1} + H_i$.  Let $x \in V_i$.  By Lemma \ref{thm:raise} and (\ref{eq:b4}), there exists $y \in V_{i-1}$ such that $[A, A^*]^{d-2i+2} y = [A, A^*]^{d-2i+1}x$.  Using this we find that $x- [A, A^*] y \in H_i$.  So, $x \in [A,A^*] V_{i-1}+H_i$.  We have now shown equality in (\ref{refineVbi2}).  It remains to show that the sum in (\ref{refineVbi2}) is  direct.  To do this we show that $[A, A^*] V_{i-1} \cap H_i=0$.  Let $x \in [A,A^*] V_{i-1} \cap H_i$.  Using (\ref{H_i1}) we find $[A,A^*]^{d-2i+1}x=0$.  Also, there exists $y \in V_{i-1}$ such that $x=[A, A^*]y$.  Combining these facts with (\ref{eq:b4}), we find $y=0$, and then $x=0$.  We have now shown that the sum in (\ref{refineVbi2}) is direct, and this completes the proof of (\ref{refineVbi2}).  Combining (\ref{refineVbi1}) and (\ref{refineVbi2}) we find
\begin{eqnarray*}
V_i=\sum_{j=0}^{i} [A, A^*]^{i-j} H_j \qquad \hbox{(direct sum)}. \\
\end{eqnarray*}
Case 2: $d/2 < i \leq d$.  This case follows immediately from Case 1 and (\ref{eq:b4}). \\ \\
We have now shown (\ref{refineV1}).  From (\ref{refineV1}) and Note \ref{A^*A} we immediately obtain (\ref{refineV2}).
\hfill $\Box$ \\

\begin{lemma}
\label{thm:Vsum}
With reference to Assumption \ref{assume}, the following holds.
\begin{eqnarray*}
V_i + \cdots + V_d= \Vb_i + \cdots + \Vb_d \qquad (0 \leq i \leq d).
\end{eqnarray*}
\end{lemma}

\noindent
{\it Proof:}  For $d \leq 1$, the result follows immediately from Definition \ref{def:bidiag}.  So assume that $d \geq 2$.  First we show that $V_i + \cdots + V_d \subseteq \Vb_i + \cdots + \Vb_d$ for $0 \leq i \leq d$.  Suppose, toward a contradiction, that there exists an $i$ such that $V_i + \cdots + V_d \nsubseteq \Vb_i + \cdots + \Vb_d$.  Define $t:=\min\{ \, i \, |\, 0 \leq i \leq d\, , \, V_i + \cdots + V_d \nsubseteq \Vb_i + \cdots + \Vb_d \,\}$ and $r:=\max\{\,i \, |\, 0 \leq i \leq d\, , \, V_i + \cdots + V_d \nsubseteq \Vb_i + \cdots + \Vb_d \,\}$.  We now show that 
\begin{eqnarray}
\label{Vsum1}
1 \leq t \leq d/2.
\end{eqnarray}
Since $V_0 + \cdots + V_d \subseteq V = \Vb_0 + \cdots + \Vb_d$, we see $t \geq 1$.  Now suppose, toward a contradiction, that $t > d/2$.  By (\ref{H_i1}) and the minimality of $t$, we have $H_k \subseteq V_k + \cdots + V_d \subseteq \Vb_k + \cdots + \Vb_d$ for $0 \leq k \leq d-t$.  This and Lemma \ref{thm:raise} give $\sum_{k=0}^{d-j} [A, A^*]^{j-k} H_k \subseteq \Vb_j + \cdots + \Vb_d$ for $t \leq j \leq d$.  From this and (\ref{refineV1}), we find $V_t + \cdots + V_d = \sum_{j=t}^d \sum_{k=0}^{d-j} [A, A^*]^{j-k} H_k  \subseteq \Vb_t + \cdots + \Vb_d$.  This contradicts the definition of $t$, and so $t \leq d/2$.  We have now shown (\ref{Vsum1}).   We now show that
\begin{eqnarray}
\label{Vsum2}
t \leq r \leq d-t.
\end{eqnarray}
It is immediate from the construction that $t \leq r$.  We now show that $r \leq d-t$.  To do this we show that 
\begin{eqnarray}
\label{Vsum3}
V_{d-t+s} + \cdots + V_d \subseteq \Vb_{d-t+s} + \cdots + \Vb_d \qquad (1 \leq s \leq t).
\end{eqnarray}
By (\ref{H_i1}) and the minimality of $t$, we find $H_k \subseteq V_k + \cdots + V_d \subseteq \Vb_k + \cdots + \Vb_d$ for $0 \leq k \leq t-1$.  From this and Lemma \ref{thm:raise} we have $\sum_{k=0}^{t-j} [A, A^*]^{d-t+j-k} H_k \subseteq \Vb_{d-t+j} + \cdots + \Vb_d$ for $s \leq j \leq t$.  This and (\ref{refineV1}) give $V_{d-t+s}+ \cdots + V_d = \sum_{j=s}^t V_{d-t+j} = \sum_{j=s}^t \sum_{k=0}^{t-j} [A, A^*]^{d-t+j-k} H_k  \subseteq \Vb_{d-t+s} + \cdots + \Vb_d$.  We have now shown (\ref{Vsum3}).  From (\ref{Vsum3}) and the maximality of $r$ we have $r \leq d-t$, and so we have shown (\ref{Vsum2}).  By the minimality of $t$, $V_t + \cdots + V_d \subseteq V_{t-1} + \cdots + V_d \subseteq \Vb_{t-1} + \cdots + \Vb_d$.  By the definition of $t$ there exists $x \in V_t + \cdots  + V_d$ with $x  \notin \Vb_t + \cdots + \Vb_d$.  Combining the previous two sentences, we find there exist $v^*_k \in V^*_k$ $(t-1 \leq k \leq d)$ such that
\begin{eqnarray}
\label{Vsum4}
x  = v^*_{t-1} + v^*_t + \cdots + v^*_d \qquad \mbox{with} \qquad v^*_{t-1} \neq 0.
\end{eqnarray}
By (\ref{Vsum2}), $r-t+2 \geq 0$.  Applying $[A, A^*]^{r-t+2}$ to (\ref{Vsum4}) leads to
\begin{eqnarray*}
[A,A^*]^{r-t+2} v^*_{t-1} = [A,A^*]^{r-t+2}x - [A,A^*]^{r-t+2} (v^*_t + \cdots + v^*_d).
\end{eqnarray*}
By Lemma \ref{thm:raise} we have $[A,A^*]^{r-t+2} (v^*_t + \cdots + v^*_d) \in \Vb_{r+2} + \cdots + \Vb_d$.  Recall $x \in V_t + \cdots + V_d$.  So the maximality of $r$ yields $[A,A^*]^{r-t+2} x \in V_{r+2} + \cdots + V_d \subseteq \Vb_{r+2} + \cdots + \Vb_d$.  Combining the previous three sentences, we have $[A,A^*]^{r-t+2} v^*_{t-1} \in \Vb_{r+2} + \cdots + \Vb_d$.  However, by Lemma \ref{thm:raise}, $[A,A^*]^{r-t+2} v^*_{t-1} \in \Vb_{r+1}$, and so  $[A,A^*]^{r-t+2} v^*_{t-1} \in \Vb_{r+1} \cap (\Vb_{r+2} + \cdots + \Vb_d)$.  From this, and since $\sum_{i=0}^d \Vb_i$ is a direct sum, we find
\begin{eqnarray}
\label{Vsum5}
[A,A^*]^{r-t+2} v^*_{t-1} =0.
\end{eqnarray}
By (\ref{Vsum2}), $d-t-r \geq 0$.  Applying $[A,A^*]^{d-t-r}$ to (\ref{Vsum5}) yields $[A,A^*]^{d-2t+2} v^*_{t-1} =0$.  This and (\ref{eq:b3}) give $v^*_{t-1} =0$, which contradicts (\ref{Vsum4}).  We have now shown that $V_i + \cdots + V_d \subseteq \Vb_i + \cdots + \Vb_d$ for $0 \leq i \leq d$.  This and Note \ref{A^*A} yield $V^*_i + \cdots + V^*_d \subseteq V_i + \cdots + V_d$ for $0 \leq i \leq d$, and we obtain the desired result.
\hfill $\Box$ \\

\noindent
We are now ready to prove Lemma \ref{thm:olddefn}. \\

\noindent
{\it Proof of Lemma \ref{thm:olddefn}:}  First we prove (\ref{olddefn1}).  Recall for $0 \leq j \leq d$, $V_j$ is the eigenspace of $A$ corresponding to eigenvalue $\theta_j$.  From this and Lemma \ref{thm:Vsum}, we have
\begin{align*}
(A - \theta_i I) \Vb_i &\subseteq (A - \theta_i I) (\Vb_i + \cdots + \Vb_d) \\
&= (A - \theta_i I) (V_i + \cdots + V_d) \\
&\subseteq V_{i+1} + \cdots + V_d \\
&= \Vb_{i+1} + \cdots + \Vb_d.
\end{align*}
By (\ref{eq:b1}) we see that $(A - \theta_i I)  \Vb_i \subseteq \Vb_i + \Vb_{i+1}$.  Combining the previous two sentences with the fact that $\sum^d_{i=0} \Vb_i$ is a direct sum, we obtain (\ref{olddefn1}).  From (\ref{olddefn1}) and Note \ref{A^*A} we immediately obtain (\ref{olddefn2}).
 \hfill $\Box$ \\

\begin{lemma}
\label{brackpoly}
With reference to Assumption \ref{assume},  the following {\rm (i)--(iv)} hold.
\begin{enumerate}
\item[\rm (i)] $[A, A^*] |_{\Vb_{i}} = (\theta^*_i - \theta^*_{i+1})(A - \theta_i I)|_{\Vb_{i}}$  \qquad $(0 \leq i \leq d)$,
\item[\rm (ii)] $[A, A^*] |_{V_{i}} = (\theta_{i+1} - \theta_{i})(\Ab - \Tb_i I)|_{V_{i}}$ \qquad $(0 \leq i \leq d)$,
\item[\rm (iii)] $[A, A^*]^k |_{\Vb_{i}} = \prod_{s=i}^{k+i-1} (\theta^*_s - \theta^*_{s+1})(A - \theta_s I)|_{\Vb_{i}}$ \qquad $(0 \leq i \leq d)$, \, $(0 \leq k \leq d-i+1)$,
\item[\rm (iv)] $[A, A^*]^k |_{V_{i}} = \prod_{s=i}^{k+i-1} (\theta_{s+1} - \theta_s)(\Ab - \Tb_s I)|_{V_{i}}$ \qquad $(0 \leq i \leq d)$, \, $(0 \leq k \leq d-i+1)$.
\end{enumerate}
\end{lemma}

\noindent
{\it Proof:} (i) Let $v \in \Vb_i$.  By (\ref{olddefn1}) there exists $v' \in \Vb_{i+1}$ such that $Av = \theta_i v + v'$.  Using this we have
\begin{align*}
[A,A^*]v &= (\theta^*_i  - \Ab)Av \\
&= \theta_i (\theta^*_i  - \Ab)v + (\theta^*_i  - \Ab)v' \\
&= (\theta^*_i  - \Tb_{i+1}) v' \\
&=  (\theta^*_i  - \Tb_{i+1}) (A - \theta_i I) v.
\end{align*}
(ii) Similar to (i). \\
(iii), (iv) Immediate from (i), (ii), and Lemma \ref{thm:olddefn}.
\hfill $\Box$ \\

\section{The projections $E_i, \, \Eb_i$}

\noindent
In this section we introduce linear transformations $E_i, \, \Eb_i$ which will be the main tools used in proving Theorem \ref{thm:main1}.

\begin{definition}
\label{def:EiEbi}
\rm
With reference to Assumption \ref{assume}, define the following linear transformations.
\begin{enumerate}
\item[\rm (i)] For $0 \leq i \leq d$, let $E_i : V \rightarrow V$ denote the linear transformation satisfying both
\begin{center}
$(E_i - I)V_i=0$, \\
$E_iV_j=0 \,\,\,\,  \mbox{for} \,\,\,\, j \neq i \,\,\,\, \mbox{and} \,\,\,\, 0 \leq j \leq d$.
\end{center}
\item[\rm (ii)] For $0 \leq i \leq d$, let $\Eb_i : V \rightarrow V$ denote the linear transformation satisfying both
\begin{center}
$(\Eb_i - I) \Vb_i=0$, \\
$\Eb_i \Vb_j=0 \,\,\,\,  \mbox{for} \,\,\,\, j \neq i \,\,\,\, \mbox{and} \,\,\,\, 0 \leq j \leq d$.
\end{center}
\end{enumerate}
In other words, $E_i$ (resp.~$\Eb_i$) is the projection map from $V$ onto $V_i$ (resp.~$\Vb_i$).  For notational convenience, let $E_{d+1}:=0$ and $\Eb_{d+1}:=0$.
\end{definition}

\begin{lemma}
\label{thm:Eprops}
With reference to Assumption \ref{assume} and Definition \ref{def:EiEbi}, the following hold.
\begin{eqnarray}
\label{Eprops1}
A \, E_i = E_i \, A = \theta_i \, E_i, & & \Ab \, \Eb_i = \Eb_i \, \Ab = \Tb_i \, \Eb_i \qquad (0 \leq i \leq d), \\
\label{Eprops2}
E_i \, E_j = \delta_{ij} \, E_i, & & \Eb_i \, \Eb_j = \delta_{ij} \, \Eb_i \qquad (0 \leq i,j \leq d), \\
\label{Eprops3}
\sum_{i=0}^d E_i=I, & & \sum_{i=0}^d \Eb_i=I, \\
\label{Eprops4}
A = \sum_{i=0}^d \theta_i E_i, & & \Ab = \sum_{i=0}^d \Tb_i \Eb_i.
\end{eqnarray}
\end{lemma}

\noindent
{\it Proof:}  Immediate from Definition \ref{def:EiEbi}.
\hfill $\Box$ \\

\begin{lemma}
\label{thm:EAEb}
With reference to Assumption \ref{assume} and Definition \ref{def:EiEbi}, the following {\rm (i)},{\rm (ii)} hold for $0 \leq i,j \leq d$.
\begin{enumerate}
\item[\rm (i)] $\Eb_j \, A \, \Eb_i = \left\{ \begin{array}{cl}
 \theta_i \Eb_i, & \mbox{if } j=i; \\ 0, &  \mbox{if } j-i > 1 \, \mbox{or } \, i-j \geq 1; \\ \neq 0, & \mbox{if } j-i=1, \end{array} \right.$
\item[\rm (ii)] $E_j \, \Ab \, E_i = \left\{ \begin{array}{cl}
 \Tb_i E_i, & \mbox{if } j=i; \\ 0, &  \mbox{if } j-i > 1 \, \mbox{or } \, i-j \geq 1; \\ \neq 0, & \mbox{if } j-i=1. \end{array} \right.$
\end{enumerate}
\end{lemma}

\noindent
{\it Proof:}  (i) Let $v \in V$, and observe $\Eb_i v \in \Vb_i$.  By (\ref{olddefn1}) we find $(A - \theta_i I) \Eb_i v = v'$, for some $v' \in \Vb_{i+1}$.  Multiplying this equation on the left by $\Eb_j$, we have
\begin{eqnarray}
\label{EAEb1}
(\Eb_j A \Eb_i - \theta_i \Eb_j \Eb_i) v = \Eb_j v'.
\end{eqnarray}
Suppose that $j=i$.  Then by (\ref{Eprops2}), (\ref{EAEb1}) becomes $(\Eb_i A \Eb_i - \theta_i \Eb_i) v = 0$.  So $\Eb_i A \Eb_i = \theta_i \Eb_i$.  Now suppose that $j-i > 1$ or $i-j \geq 1$.  Then by (\ref{Eprops2}), (\ref{EAEb1}) becomes $\Eb_j A \Eb_i v = 0$.  So $\Eb_j A \Eb_i = 0$.  Now suppose that $j-i=1$.  We show that $\Eb_j \, A \, \Eb_i$ is not the zero transformation.  Combining Lemma \ref{thm:injsurj} with Lemma \ref{brackpoly}(i), we find that for $0 \leq i < d/2$ (resp.~$d/2 \leq i \leq d$), $(A - \theta_i I) |_{\Vb_i} : \Vb_i \rightarrow \Vb_{i+1}$ is an injection (resp.~surjection).  Thus, for $0 \leq i \leq d$, $(A - \theta_i I) |_{\Vb_i} : \Vb_i \rightarrow \Vb_{i+1}$ is not the zero transformation.  So there exists $w \in V^*_i$ such that $(A - \theta_i I) \Eb_i w = w'$, for some $0 \neq w' \in \Vb_{i+1}$.  Multiplying this equation on the left by $\Eb_j$ yields $\Eb_j A \Eb_i w = \Eb_j w' = w'$, since $j=i+1$.  Therefore, $\Eb_j \, A \, \Eb_i$ is not the zero transformation. \\
(ii) Immediate from (i) and Note \ref{A^*A}.
\hfill $\Box$ \\

\begin{lemma}
\label{thm:EAr}
With reference to Assumption \ref{assume} and Definition \ref{def:EiEbi}, the following {\rm(i)},{\rm(ii)} hold.
\begin{enumerate}
\item[\rm (i)] $E_j \, A^{* \, r} \, E_i = 0$ \qquad $(0 \leq i,j \leq d)$, \, $(0 \leq r < j-i)$,
\item[\rm (ii)] $E_j \, A^{* \, r} \, E_i \neq 0$ \qquad $(0 \leq i,j \leq d)$, \, $(0 \leq r = j-i)$.
\end{enumerate}
\end{lemma}

\noindent
{\it Proof:}  (i) By (\ref{eq:b2}) we have
\begin{eqnarray}
\label{EAr1}
A^{* \, r} V_i \subseteq V_i + V_{i+1} + \cdots + V_{i+r}.
\end{eqnarray}
By Definition \ref{def:EiEbi}(i), and since $r+i < j$, we have
\begin{eqnarray}
\label{EAr2}
E_j V_i + E_j V_{i+1} + \cdots + E_j V_{i+r} = 0.
\end{eqnarray}
Let $v \in V$, and observe $E_i v \in V_i$.  Combining (\ref{EAr1}) and (\ref{EAr2}) gives $E_j A^{* \, r} E_i v = 0$.  So $E_j A^{* \, r} E_i = 0$. \\
(ii) Combining Lemma \ref{brackpoly}(iv) with (i), and since $r=j-i$, we have
\begin{eqnarray}
\label{EAr3}
E_j \, [A,A^*]^r \, E_i = \Big( \prod_{s=i}^{r+i-1} (\Tb_{s+1} - \Tb_s) \Big) \,\, E_j A^{* \, r} E_i.
\end{eqnarray}
We now divide the proof into two cases.  \\ \\
Case 1:  $0 \leq j < d-i$.  Since $j=i+r$, we see that $0 \leq i < d/2$ and $0 \leq r < d-2i$.  Thus, by Lemma \ref{thm:injsurj}, we find that the restriction $[A,A^*]^r |_{V_i} : V_i \rightarrow V_j$ is an injection.  Let $0 \neq v \in V_i$.  Then $E_j [A,A^*]^r E_i v \neq 0$.  Combining this with (\ref{EAr3}) yields $E_j A^{* \, r} E_iv \neq 0$.  So $E_j A^{* \, r} E_i$ is not the zero transformation. \\ \\
Case 2:  $d-i \leq j \leq d$.  First we show that
\begin{eqnarray}
\label{EAr4}
[A,A^*]^r |_{V_i} : V_i \rightarrow V_j \,\,\, {\rm is \,\, a \,\, surjection}.
\end{eqnarray}
Suppose that $0 \leq i \leq d/2$.  By (\ref{eq:b4}), we have $[A,A^*]^{d-2i} |_{V_i} : V_i \rightarrow V_{d-i}$ is a bijection.  By Lemma \ref{thm:injsurj}, we have $[A,A^*]^{j+i-d} |_{V_{d-i}} : V_{d-i} \rightarrow V_j$ is a surjection.  Combining the previous two sentences we obtain (\ref{EAr4}).  Now suppose that $d/2 < i \leq d$.  Then (\ref{EAr4}) follows immediately from Lemma \ref{thm:injsurj}.  We have now shown (\ref{EAr4}).  Let $0 \neq v \in V_j$.  Then by (\ref{EAr4}), there exists $v' \in V_i$ such that $E_j [A,A^*]^r E_i v' = E_j v = v \neq 0$.  Combining this with (\ref{EAr3}) yields $E_j A^{* \, r} E_iv' \neq 0$.  So $E_j A^{* \, r} E_i$ is not the zero transformation.
\hfill $\Box$ \\

\section{The proof of the fundamental bidiagonal relation}

\noindent
In this section we prove Theorem \ref{thm:main1}.  The following two lemmas will be used in proving Theorem \ref{thm:main1}.

\begin{lemma}
\label{thm:relrec}
Let $b, \, \alpha, \, \ab, \, \gamma$ denote scalars in $\K$.  With reference to Assumption \ref{assume}, suppose that $d \geq 1$.  Then
\begin{eqnarray}
\label{relrec1}
A \Ab - b \Ab A - \alpha A - \ab \Ab - \gamma I & = & 0
\end{eqnarray}
if and only if the following {\rm (i)--(iii)} below hold.
\begin{enumerate}
\item[\rm (i)] $\theta_{i+1} - b \theta_i = \ab$ \qquad $(0 \leq i \leq d-1)$,
\item[\rm (ii)] $\Tb_i - b \Tb_{i+1} = \alpha$ \qquad $(0 \leq i \leq d-1)$,
\item[\rm (iii)] $b \theta_i \Tb_{i+1} - \theta_{i+1} \Tb_i = \gamma$ \qquad $(0 \leq i \leq d-1)$.
\end{enumerate}
\end{lemma}

\noindent
{\it Proof:}  Let $E_i, \, \Eb_i$ be as in Definition \ref{def:EiEbi}.  Let $C$ denote the expression on the left in (\ref{relrec1}). For $0 \leq i,j \leq d$, we evaluate $E_j C E_i$ using (\ref{Eprops1}) and get
\begin{eqnarray}
\label{relrec2}
E_jCE_i = (\theta_j - b \theta_i - \ab) E_j \Ab E_i - (\alpha \theta_j + \gamma) E_jE_i.
\end{eqnarray}
For $0 \leq i,j \leq d$, we evaluate $\Eb_j C \Eb_i$ using (\ref{Eprops1}) and get
\begin{eqnarray}
\label{relrec3}
\Eb_j C \Eb_i = (\Tb_i - b \Tb_j - \alpha) \Eb_j A \Eb_i - (\ab \Tb_j + \gamma) \Eb_j \Eb_i.
\end{eqnarray}
$(\Longrightarrow)$:  Assume (\ref{relrec1}) holds, so that $C=0$.  We show that (i)--(iii) hold.  We have $E_{i+1} C E_i = 0$ since $C=0$.  So by (\ref{relrec2}) and (\ref{Eprops2}), we find $(\theta_{i+1} - b \theta_i - \ab) E_{i+1} \Ab E_i = 0$.  By Lemma \ref{thm:EAEb}(ii) we have $E_{i+1} \Ab E_i \neq 0$, and so $\theta_{i+1} - b \theta_i - \ab = 0$.  We have now shown (i).  We have $\Eb_{i+1} C \Eb_i = 0$ since $C=0$.  So, by (\ref{relrec3}) and (\ref{Eprops2}), we find $(\Tb_i - b \Tb_{i+1} - \alpha) \Eb_{i+1} A \Eb_i = 0$.  By Lemma \ref{thm:EAEb}(i) we have $\Eb_{i+1} A \Eb_i \neq 0$, and so $\Tb_i - b \Tb_{i+1} - \alpha = 0$.  We have now shown (ii).  We know $E_i C E_i = 0$ since $C=0$.  So by (\ref{relrec2}), (\ref{Eprops2}), and Lemma \ref{thm:EAEb}(ii), we have $( (\theta_i - b \theta_i - \ab) \Tb_i - \alpha \theta_i - \gamma) E_i = 0$.  From this and since $E_i \neq 0$, we find $ (\theta_i - b \theta_i - \ab) \Tb_i - \alpha \theta_i = \gamma$.  Substituting (i),(ii) into this equation, and simplifying gives $b \theta_i \Tb_{i+1} - \theta_{i+1} \Tb_i = \gamma$.  We have now shown (iii).  We have now shown (i)--(iii) hold. \\ \\
$(\Longleftarrow)$:  Assume (i)--(iii) hold.  We show that $C=0$.  By (\ref{Eprops3}) we have $C = \sum_{j=0}^d \sum_{i=0}^d E_jCE_i$.  So to show that $C=0$, it suffices to show $E_jCE_i=0$ for $0 \leq i,j\leq d$.  We divide the argument into three cases. \\
Case 1:  $j-i > 1$ or $i-j \geq 1$.  Simplifying (\ref{relrec2}), using (\ref{Eprops2}) and Lemma \ref{thm:EAEb}(ii), gives $E_jCE_i=0$.  \\ \\
Case 2:  $j-i = 1$.  Simplifying (\ref{relrec2}), using (\ref{Eprops2}) and (i), gives $E_jCE_i=0$. \\ \\
Case 3:  $j = i$.  Observe (iii) is equivalent to $\theta_i \Tb_i - b \theta_i \Tb_i - (\theta_{i+1} - b \theta_i) \Tb_i -(\Tb_i - b \Tb_{i+1}) \theta_i = \gamma$.  Substituting (i),(ii) into this equation yields $(\theta_i - b \theta_i - \ab) \Tb_i - \alpha \theta_i - \gamma = 0$.  Simplifying (\ref{relrec2}), using this, (\ref{Eprops2}), and Lemma \ref{thm:EAEb}(ii), gives $E_jCE_i=0$.  We have now shown that $E_jCE_i=0$ for $0 \leq i,j\leq d$.  So $C=0$, and (\ref{relrec1}) holds.
\hfill $\Box$ \\

\begin{lemma}
\label{thm:ghpoly}
With reference to Assumption \ref{assume}, suppose that $d \geq 2$.  Then there exist polynomials $g,\, h \in \K[\lambda]$ with $1 \leq \deg \, g \leq d-1$ and $0 \leq \deg \, h \leq d$ such that
\begin{eqnarray}
\label{ghpoly1}
A \Ab - g(\Ab) A - h(\Ab) = 0.
\end{eqnarray}
\end{lemma}

\noindent
{\it Proof:}  Let $\lbrace \theta_i \rbrace _{i=0}^d$ and $\lbrace \Tb_i \rbrace _{i=0}^d$ be as in Assumption \ref{assume}.  First we construct the polynomial $g$.  Define the scalars $M_{i,j} := \theta^{* \, j-1}_i$ for $1 \leq i,j \leq d$.  Let $M$ denote the $d$ by $d$ matrix with $(i,j)$ entry $M_{i,j}$ for $1 \leq i,j \leq d$.  Since $\Tb_0, \, \Tb_1, \, \ldots, \, \Tb_d$ are distinct, and $M$ is a Vandermonde matrix, we find that $M$ is invertible.  For $0 \leq i \leq d-1$, define the scalars $g_i$ by the following matrix equation:  $[g_0, g_1, \ldots g_{d-1}]^T := M^{-1} [\Tb_0, \Tb_1, \ldots, \Tb_{d-1}]^T$.  Define $g \in \K[\lambda]$ as $g := \sum_{s=0}^{d-1} g_s \lambda ^s$, and observe $\deg g \leq d-1$.  By construction $M [g_0, g_1, \ldots g_{d-1}]^T = [\Tb_0, \Tb_1, \ldots, \Tb_{d-1}]^T$, and so 
\begin{eqnarray}
\label{ghpoly2}
g(\Tb_{i+1}) = \Tb_i \qquad (0 \leq i \leq d-1).
\end{eqnarray}
By (\ref{ghpoly2}) and since $\Tb_0, \, \Tb_1, \, \ldots, \, \Tb_d$ are distinct, we see that $g$ is not a constant polynomial.  So $\deg g \geq 1$.  Now we construct the polynomial $h$.  Define the scalars $N_{i,j} := \theta^{* \, j}_i$ for $0 \leq i,j \leq d$.  Let $N$ denote the $(d+1)$ by $(d+1)$ matrix with $(i,j)$ entry $N_{i,j}$ for $0 \leq i,j \leq d$.  Since $\Tb_0, \, \Tb_1, \, \ldots, \, \Tb_d$ are distinct, and $N$ is a Vandermonde matrix, we find that $N$ is invertible.  For $0 \leq i \leq d$, define the scalars $h_i$ by the following matrix equation: \\ $[h_0, h_1, \ldots h_d]^T := N^{-1} [\theta_0(\Tb_0 - g(\Tb_0)), \, \theta_1(\Tb_1 - g(\Tb_1)), \, \ldots, \, \theta_d(\Tb_d - g(\Tb_d))]^T$.  Define $h \in \K[\lambda]$ as $h := \sum_{s=0}^{d} h_s \lambda ^s$, and observe $0 \leq \deg h \leq d$.  By construction 
\begin{eqnarray}
\label{ghpoly3}
h(\Tb_i) = \theta_i (\Tb_i - g(\Tb_i)) \qquad (0 \leq i \leq d).
\end{eqnarray}
We now show (\ref{ghpoly1}).  Let $v \in \Vb_i$, and recall that $v$ is an eigenvector for $\Ab$ with eigenvalue $\Tb_i$.  Also, recall $(A - \theta_i I)v \in \Vb_{i+1}$ by (\ref{olddefn1}).  Using this, $(\ref{ghpoly2})$, and $(\ref{ghpoly3})$, we have
\begin{eqnarray*}
(A \Ab - g(\Ab)A - h(\Ab) )v
&=& ( (A - \theta_i I) \Ab + \theta_i \Ab - g(\Ab)(A- \theta_i I) - \theta_i g(\Ab) - h(\Ab) )v \\
&=& ( \Tb_i (A-\theta_i I) + \Tb_i \theta_i - g(\Tb_{i+1})(A-\theta_i I) - \theta_i g(\Tb_i) - h(\Tb_i) )v \\ &=& 0.
\end{eqnarray*}
We have now shown that the left-hand side of (\ref{ghpoly1}) vanishes on $\Vb_i$ for $0 \leq i \leq d$.  From this we obtain (\ref{ghpoly1}), since $\lbrace V^*_i \rbrace _{i=0}^d$ is a decomposition of $V$.
\hfill $\Box$ \\

\noindent
We are now ready to prove Theorem \ref{thm:main1}. \\

\noindent
{\it Proof of Theorem \ref{thm:main1}:}
First assume that $d \geq 2$.  By Lemma \ref{thm:ghpoly} there exist polynomials $g,\, h \in \K[\lambda]$ with $1 \leq \deg g \leq d-1$ and $0 \leq \deg h \leq d$ such that
\begin{eqnarray}
\label{main11}
A \Ab - g(\Ab) A - h(\Ab) =0.
\end{eqnarray}
Let $k := \max (\deg g, \, \deg h)$, and observe that $1 \leq k \leq d$.  We now show
\begin{eqnarray}
\label{main11.5}
k = 1.
\end{eqnarray}
Suppose, towards a contradiction, that $k > 1$.  Let $g_k$ (resp.~$h_k$) denote the coefficient of $\lambda ^k$ in $g$ (resp.~$h$).  For $0 \leq i \leq d$, let $E_i, \, \Eb_i$ be as in Definition \ref{def:EiEbi}.  Multiplying each term in (\ref{main11}) on the left by $E_k$, and on the right by $E_0$, and evaluating the result using (\ref{Eprops1}) and Lemma \ref{thm:EAr}(i) yields $0 = (\theta_0 g_k + h_k )E_k A ^{* \, k} E_0$.  From this and Lemma \ref{thm:EAr}(ii), we have
\begin{eqnarray}
\label{main12}
0 = \theta_0 g_k + h_k.
\end{eqnarray}
We now break the argument into two cases. \\ \\
Case 1:  $k=d$.  Since $\deg g \leq d-1$, we see $g_k=0$.  From this and (\ref{main12}), we see $h_k=0$. Combining the previous two sentences, we find $k = \max (\deg g, \, \deg h) \leq k-1$, for a contradiction.  So (\ref{main11.5}) holds.  \\ \\
Case 2:  $k \leq d-1$.  Multiplying each term in (\ref{main11}) on the left by $E_{k+1}$, and on the right by $E_1$, and evaluating the result using (\ref{Eprops1}) and Lemma \ref{thm:EAr}(i), we have
$0 = (\theta_1 g_k + h_k )E_{k+1} A ^{* \, k} E_1$.  From this and Lemma \ref{thm:EAr}(ii), we have
\begin{eqnarray}
\label{main13}
0 = \theta_1 g_k + h_k.
\end{eqnarray}
Combining (\ref{main12}) and (\ref{main13}) yields $(\theta_0 - \theta_1) g_k = 0$.  Thus, since $\theta_0 \neq \theta_1 $, we see $g_k = 0$.  From this and (\ref{main12}), we see $h_k =0$.  Combining the previous two sentences, we find $k = \max (\deg g, \, \deg h) \leq k-1$, for a contradiction.  So (\ref{main11.5}) holds. \\ \\ We have now shown (\ref{main11.5}), and so $\deg g = 1$ and $\deg h \leq 1$.  Therefore, there exist scalars $b, \, \alpha, \, \ab, \, \gamma$ in $\K$ with $b$ nonzero such that $g = b \lambda + \alpha$ and $h = \ab \lambda + \gamma$.  From this and (\ref{main11}), we obtain (\ref{main1.1}).  We now show that the sequence of scalars $b, \, \alpha, \, \ab, \, \gamma$ is uniquely determined by the pair $A, \, \Ab$.  Let $b, \, \alpha, \, \ab, \, \gamma$ denote any sequence of scalars in $\K$ which satisfies (\ref{main1.1}).  Applying Lemma \ref{thm:relrec}, $\theta_{i+1} - b \theta_i = \ab$ for $0 \leq i \leq d-1$.  This gives $b = (\theta_{i+1} - \theta_i) / (\theta_i - \theta_{i-1})$ for $1 \leq i \leq d-1$.  So $b$ is uniquely determined by $A, \, A^*$.  From this and Lemma \ref{thm:relrec}(i)--(iii), we find $\ab, \, \alpha, \, \gamma$ are also uniquely determined by $A, \, A^*$.  We have now proved the theorem for the case $d \geq 2$.  \\ \\
Now assume that $d=1$.  Let $b$ denote any nonzero scalar in $\K$.  Define $\ab := \theta_1 - b \theta_0$, $\alpha := \Tb_0 - b \Tb_1$, and $ \gamma := b \theta_0 \Tb_1 - \theta_1 \Tb_0$.  Applying Lemma \ref{thm:relrec}, we find $b, \, \alpha, \, \ab, \, \gamma$ satisfy
(\ref{main1.1}).  \\ \\
Now assume that $d=0$.  Let $b$ denote any nonzero scalar in $\K$, and let $\alpha, \, \ab$ denote any scalars in $\K$.  Define $\gamma := (\theta_0 - b \theta_0 - \ab) \Tb_0 - \alpha \theta_0$.  Let $C$ denote the left-hand side of (\ref{main1.1}).  By (\ref{Eprops3}) we have $C = E_0 C E_0$.  Evaluating $E_0 C E_0$ using (\ref{Eprops1}), (\ref{Eprops2}), and Lemma \ref{thm:EAEb}(ii) yields $C = ((\theta_0 - b \theta_0 - \ab) \Tb_0 - \alpha \theta_0 - \gamma) E_0$.  From this, and the definition of $\gamma$, we see that $C=0$.  Hence, (\ref{main1.1}) holds.
\hfill $\Box$ \\

\section{The proof of the reducibility of bidiagonal pairs}

\noindent
In this section we prove Theorem \ref{thm:main2}.  The following two lemmas will be used in the proof of Theorem \ref{thm:main2}, and in the proof of Theorem \ref{thm:class}.

\begin{lemma}
\label{thm:main2prep}
Adopt Assumption \ref{assume}.  Suppose that $d \geq 2$, and let $b$ denote the base of $A, \, A^*$.  Then the following holds.
\begin{eqnarray*}
\frac{\theta_{i+1} - \theta_i}{\theta_i - \theta_{i-1}} = b = \frac{\Tb_i - \Tb_{i-1}}{\Tb_{i+1} - \Tb_i} \qquad (1 \leq i \leq d-1).
\end{eqnarray*}
\end{lemma}

\noindent
{\it Proof:}  Let $b, \, \alpha, \, \alpha^*$ be as in (\ref{main1.1}).  Combining Theorem \ref{thm:main1} and Lemma \ref{thm:relrec}, we have $\theta_{i+1} - b \theta_i = \ab$ and $\Tb_i - b \Tb_{i+1} = \alpha$ for $0 \leq i \leq d-1$.  From this, we find $\theta_{i+1} - b \theta_i = \theta_{i} - b \theta_{i-1}$ and $\Tb_{i} - b \Tb_{i+1} = \Tb_{i-1} - b \Tb_{i}$ for $1 \leq i \leq d-1$.  Solving these two equations for $b$ we obtain the desired result.
\hfill $\Box$ \\

\begin{lemma}
\label{thm:main2prep2}
Let $d$ denote a nonnegative integer.  Let $\{ \beta_i \}_{i=0}^d$ and $\{ \beta^*_i \}_{i=0}^d$ denote two sequences of scalars in $\K$ such that $\beta_i \neq \beta_j$ and $\beta^*_i \neq \beta^*_j$ for $0 \leq i,j \leq d$ with $i \neq j$.  For $d \leq 1$, define $b := 1$.  For $d \geq 2$, assume that the expressions
\begin{eqnarray*}
\frac{\beta_{i+1} - \beta_i}{\beta_i - \beta_{i-1}}, \qquad \frac{\beta^*_i - \beta^*_{i-1}}{\beta^*_{i+1} - \beta^*_i}
\end{eqnarray*}
are equal and independent of $i$ for $1 \leq i \leq d-1$, and define $b$ as their common value.  Suppose that $b = 1$.  Then there exists scalars $b_1, \, b_2, \, c_1, \, c_2$ in $\K$ with $b_2, \, c_2$ both nonzero such that 
\begin{eqnarray}
\label{main2.1prep2}
\beta_i &=& b_1 + b_2 2i \qquad (0 \leq i \leq d), \\
\label{main2.2prep2}
\beta_i ^* &=& c_1 + c_2 (-2i) \qquad (0 \leq i \leq d).
\end{eqnarray}
Suppose that $b \neq 1$.  Let $q$ denote a scalar in $\K$ such that $b=q^{-2}$. Then there exists scalars $b_1 ', \, b_2 ', \, c_1 ', \, c_2 '$ in $\K$ with $b_2 ', \, c_2 '$ both nonzero such that 
\begin{eqnarray}
\label{main2.3prep2}
\beta_i &=& b_1 '+ b_2 ' \, q^{-2i} \qquad (0 \leq i \leq d), \\
\label{main2.4prep2}
\beta^*_i &=& c_1 ' + c_2 ' \, q^{2i} \qquad (0 \leq i \leq d).
\end{eqnarray}
\end{lemma}

\noindent
{\it Proof:}  First assume that $d \geq 2$.  Suppose that $b = 1$.  Solving the two recurrence relations
\begin{eqnarray*}
\frac{\beta_{i+1} - \beta_i}{\beta_i - \beta_{i-1}} = 1 = \frac{\beta^*_i - \beta^*_{i-1}}{\beta^*_{i+1} - \beta^*_i}  \qquad (1 \leq i \leq d-1),
\end{eqnarray*}
we find that there exists scalars $b_1, \, b_2, \, c_1, \, c_2$ in $\K$ with $b_2, \, c_2$ nonzero satisfying (\ref{main2.1prep2}), (\ref{main2.2prep2}).  Now suppose that $b \neq 1$.  Solving the two recurrence relations
\begin{eqnarray*}
\frac{\beta_{i+1} - \beta_i}{\beta_i - \beta_{i-1}} = q^{-2} = \frac{\beta^*_i - \beta^*_{i-1}}{\beta^*_{i+1} - \beta^*_i} \qquad (1 \leq i \leq d-1),
\end{eqnarray*}
we find that there exists scalars $b_1 ', \, b_2 ', \, c_1 ', \, c_2 '$ in $\K$ with $b_2 ', \, c_2 '$ nonzero satisfying (\ref{main2.3prep2}), (\ref{main2.4prep2}).  We have now proved the Lemma for the case $d \geq 2$.  \\ \\
Now assume $d=1$, so that $b =1$.  Define $b_1 := \beta_0$, $b_2 := (\beta_1 - \beta_0)/2$, $c_1 := \beta^*_0$, $c_2 := (\beta^*_0 - \beta^*_1)/2$.  Since $\beta_0 \neq \beta_1$ and $\beta^*_0 \neq \beta^*_1$, we have that $b_2, \, c_2$ are both nonzero.  From these definitions we obtain (\ref{main2.1prep2}), (\ref{main2.2prep2}).  \\Ê\\
Now assume $d=0$, so that $b =1$.  Let $b_2, \, c_2$ denote any nonzero scalars in $\K$, and define $b_1 := \beta_0$, $c_1 := \beta^*_0$.  From these definitions we obtain (\ref{main2.1prep2}), (\ref{main2.2prep2}).
\hfill $\Box$ \\

\noindent
We are now ready to prove Theorem \ref{thm:main2}. \\

\noindent
{\it Proof of Theorem \ref{thm:main2}:}  Adopt Assumption \ref{assume}, and let $b$ denote the base of $A, \, A^*$.  \\ \\
First assume that $b=1$.  Combining Lemma \ref{thm:main2prep} and Lemma \ref{thm:main2prep2}, we find that there exist scalars $b_1, \, b_2, \, c_1, \, c_2$ in $\K$ with $b_2, \, c_2$ both nonzero such that 
\begin{eqnarray}
\label{2.1}
\theta_i &=& b_1 + b_2 2i \qquad (0 \leq i \leq d), \\
\label{2.2}
\theta_i ^* &=& c_1 + c_2 (-2i) \qquad (0 \leq i \leq d).
\end{eqnarray}
Define the polynomials $\sigma, \, \tau \in \K[\lambda]$ as follows:
\begin{eqnarray}
\label{sigma1}
\sigma &:=& b_2^{-1} \lambda - b_1 b_2^{-1} - d, \\
\label{tau1}
\tau &:=& c_2^{-1} \lambda - c_1 c_2^{-1} + d.
\end{eqnarray}
By Lemma \ref{affineisbidiag}, $\sigma(A), \, \tau(\Ab)$ is a bidiagonal pair on $V$ with eigenvalue sequence $\lbrace \sigma (\theta_i) \rbrace _{i=0}^d$ and dual eigenvalue sequence $\lbrace \tau (\Tb_i) \rbrace _{i=0}^d$.  Using (\ref{sigma1}) and (\ref{tau1}), we have
\begin{eqnarray}
\label{sigmaA}
A &=& b_2 \sigma(A) + (b_1 + b_2 d)I, \\
\label{tauA}
A^* &=& c_2 \tau(A^*) + (c_1 - c_2 d)I.
\end{eqnarray}
Recall $b_2, \, c_2$ are both nonzero.  From (\ref{sigmaA}), (\ref{tauA}) we find that $A, \, A^*$ is affine equivalent to $\sigma(A), \, \tau(A^*)$ (taking $\mu$ from Definition \ref{def:isom} to be the identity map on $V$). Substituting (\ref{sigmaA}) and (\ref{tauA}) into (\ref{main1.1}), and simplifying we find that $\sigma(A), \, \tau(A^*)$ has base $1$.  From (\ref{2.1})--(\ref{tau1}), we have $\sigma (\theta_i) = 2i-d$ and $\tau (\Tb_i) = d-2i$ for $0 \leq i \leq d$.  Thus, $\sigma(A), \, \tau(A^*)$ is reduced.  We have now shown that for $b=1$, $A, \, A^*$ is affine equivalent to a reduced bidiagonal pair.  \\ \\
Now assume that $b \neq 1$.  Combining Lemma \ref{thm:main2prep} and Lemma \ref{thm:main2prep2}, we find that there exist scalars $b_1 ', \, b_2 ', \, c_1 ', \, c_2 '$ in $\K$ with $b_2 ', \, c_2 '$ both nonzero such that
\begin{eqnarray}
\label{2.3}
\theta_i &=& b_1 ' + b_2 ' \, q^{-2i} \qquad (0 \leq i \leq d), \\
\label{2.4}
\theta_i ^* &=& c_1 ' + c_2 ' \, q^{2i} \qquad (0 \leq i \leq d).
\end{eqnarray}
Define the polynomials $\sigma', \, \tau' \in \K[\lambda]$ as follows:
\begin{eqnarray}
\label{sigma2}
\sigma ' &:=& b_2'^{-1} q^{d} (\lambda - b_1 '), \\
\label{tau2}
\tau ' &:=& c_2'^{-1} q^{-d} (\lambda - c_1 ').
\end{eqnarray}
Observe $b_2'^{-1} q^{d}, \, c_2'^{-1} q^{-d}$ are both nonzero.  So by Lemma \ref{affineisbidiag}, $\sigma'(A), \, \tau'(\Ab)$ is a bidiagonal pair on $V$ with eigenvalue sequence $\lbrace \sigma' (\theta_i) \rbrace _{i=0}^d$ and dual eigenvalue sequence $\lbrace \tau' (\Tb_i) \rbrace _{i=0}^d$.  Using (\ref{sigma2}) and (\ref{tau2}), we have
\begin{eqnarray}
\label{sigmaprimeA}
A &=& q^{-d} b_2' \sigma'(A) + b_1'I, \\
\label{tauprimeA}
A^* &=& q^{d} c_2' \tau'(A^*) + c_1 ' I.
\end{eqnarray}
Observe $q^{-d} b_2', \, q^d c_2'$ are both nonzero.  From (\ref{sigmaprimeA}), (\ref{tauprimeA}) we find that $A, \, A^*$ is affine equivalent to $\sigma'(A), \, \tau'(A^*)$ (taking $\mu$ from Definition \ref{def:isom} to be the identity map on $V$). Substituting (\ref{sigmaprimeA}) and (\ref{tauprimeA}) into (\ref{main1.1}), and simplifying we find that $\sigma'(A), \, \tau'(A^*)$ has base $q^{-2}$, which is not $1$.  Using (\ref{2.3})--(\ref{tau2}), we have $\sigma' (\theta_i) = q^{d-2i}$ and $\tau' (\Tb_i) = q^{2i-d}$ for $0 \leq i \leq d$.  Thus, $\sigma'(A), \, \tau'(A^*)$ is reduced.  We have now shown that for $b \neq 1$, $A, \, A^*$ is affine equivalent to a reduced bidiagonal pair.
\hfill $\Box$ \\

\section{The proof that equitable pairs act as bidiagonal pairs}

In this section we prove Theorem \ref{thm:uqsl2bidiag}. \\

\noindent
{\it Proof of Theorem \ref{thm:uqsl2bidiag}:}  Let $V$ denote a segregated $\SL$-module.  Let $Y, \, Z$ denote an equitable pair in $\SL$.  We show that the action of $Y, \, Z$ on $V$ is a reduced bidiagonal pair with base $1$.  By Lemma \ref{thm:compred}, $V$ is a direct sum of irreducible $\SL$-submodules.  Let $W$ denote one of the irreducible submodules in this sum.  First we show that the action of $Y, \, Z$ on $W$ is a bidiagonal pair.  By Lemma \ref{thm:sl2mods}, $W$ is isomorphic to $V(d)$ for some nonnegative $d$.  Let $\lbrace v_i \rbrace _{i=0}^d$ denote the basis for $W$  from Lemma \ref{thm:sl2mods}.  Identify the copy of $\SL$ given in Definition \ref{def:usl2} with the copy given in Theorem \ref{thm:ueq}, via the automorphism given in Theorem \ref{thm:ueq}.  Under this automorphism, $Z$ is mapped to $h$, and $Y$ is mapped to $-2f-h$.  Using this and Lemma \ref{thm:sl2mods}, we find that the action of $Z$ on $W$ is diagonalizable with eigenvalues $\lbrace d-2i \rbrace _{i=0}^d$ corresponding to eigenvectors $\lbrace v_i \rbrace _{i=0}^d$.   We also have 
\begin{eqnarray}
\label{sl2bidiag1}
Y. v_i = (2i-d)  v_i -2(i+1) v_{i+1} \qquad (0 \leq i \leq d).
\end{eqnarray}
The previous two sentences give $[Y,Z]. v_i = -4(i+1) v_{i+1}$ for $0 \leq i \leq d$.  This yields 
\begin{eqnarray}
\label{sl2bidiag1.5}
[Y,Z]^{d-2i}. v_i = (-4)^{d-2i}(i + d-2i) \cdots (i+1) v_{d-i} \qquad (0 \leq i \leq d/2).
\end{eqnarray}
Define $u_i:=v_{d-i}$ for $0 \leq i \leq d$.  Observe that $\lbrace u_i \rbrace _{i=0}^d$ is a basis for $W$.  By Lemma \ref{thm:sl2mods}, we have $h.u_i = (2i-d) u_i$ and $e.u_i = (i+1) u_{i+1}$ for $0 \leq i \leq d$.  Composing the automorphism $\psi$ from Remark \ref{def:psi} with the automorphism given in Theorem \ref{thm:ueq}, we obtain an automorphism of $\SL$ which maps $Y$ to $h$, and $Z$ to $2e-h$.  Using this automorphism, we find that the action of $Y$ on $W$ is diagonalizable with eigenvalues $\lbrace 2i-d \rbrace _{i=0}^d$ corresponding to eigenvectors $\lbrace u_i \rbrace _{i=0}^d$.  We also have
\begin{eqnarray}
\label{sl2bidiag2}
Z. u_i = (d-2i)  u_i + 2(i+1) \, u_{i+1} \qquad (0 \leq i \leq d).
\end{eqnarray}
The previous two sentences give $[Y,Z]. u_i = 4(i+1) u_{i+1}$ for $0 \leq i \leq d$.  This yields 
\begin{eqnarray}
\label{sl2bidiag2.5}
[Y,Z]^{d-2i}. u_i = 4^{d-2i}(i + d-2i) \cdots (i+1) u_{d-i} \qquad (0 \leq i \leq d/2).
\end{eqnarray}
Combining (\ref{sl2bidiag1})--(\ref{sl2bidiag2.5}) we find that the action of $Y, \, Z$ on $W$ is a bidiagonal pair.  We now show that the action of $Y, \, Z$ on $V$ is a bidiagonal pair.  Let $V = \sum_{j=0}^{r} V(d_j)$ denote the direct sum decomposition of $V$ into irreducible $\SL$-submodules.  Without loss of generality, $d_0 \geq d_1 \geq \cdots \geq d_r$.  For $0 \leq j \leq r$, our work above has shown that the action of $Y$ (resp.~$Z$) on $V(d_j)$ is diagonalizable with eigenvalues $\lbrace 2i-d_j \rbrace _{i=0}^{d_j}$ (resp.~$\lbrace d_j - 2i \rbrace _{i=0}^{d_j}$).  Since $V$ is segregated, then either $\lbrace d_j \rbrace _{j=0}^r \subseteq 2 \Z$ or $\lbrace d_j \rbrace _{j=0}^r \subseteq 2 \Z +1$.  Combining the previous two sentences, we find that the action of $Y$ (resp.~$Z$) on $V$ is diagonalizable with eigenvalues $\lbrace 2i-d_0 \rbrace _{i=0}^{d_0}$ (resp.~$\lbrace d_0 - 2i \rbrace _{i=0}^{d_0}$).  For $0 \leq i \leq d_0$, let $V_i$ (resp.~$V_i^*$) denote the eigenspace of $Y$ (resp.~$Z$) corresponding to $2i - d_0$ (resp.~$d_0 - 2i$).  By Lemma \ref{thm:sl2mods}, we have $\dim(V_i) = \dim(V_{d-i})$ and $\dim(V^*_i) = \dim(V^*_{d-i})$ for $0 \leq i \leq d$.  Since (\ref{sl2bidiag1})--(\ref{sl2bidiag2.5}) hold on each $V(d_j)$, we have $Y, \, Z, \, \lbrace V_i \rbrace _{i=0}^{d_0}, \, \lbrace \Vb_i \rbrace _{i=0}^{d_0}$ satisfy (\ref{eq:b1})--(\ref{eq:b4}).  Thus, the action of $Y, \, Z$ on $V$ is a bidiagonal pair with eigenvalue (resp.~dual eigenvalue) sequence $\lbrace 2i-d_0 \rbrace _{i=0}^{d_0}$ (resp.~$\lbrace d_0-2i \rbrace _{i=0}^{d_0}$).  Comparing (\ref{equitablepres}) and (\ref{main1.1}), we find that the base of $Y, \, Z$ is $1$.  The previous two sentences show that $Y, \, Z$ is reduced.  \\

\noindent
Now let $V$ denote a segregated $\uq$-module.  Let $y, \, z$ denote an equitable pair in $\uq$.  We show that the action of $y, \, z$ on $V$ is a reduced bidiagonal pair with base not equal to $1$.  By Lemma \ref{thm:uqcompred}, $V$ is a direct sum of irreducible $\uq$-submodules.  Let $W$ denote one of the irreducible submodules in this sum.  First we show that the action of $y, \, z$ on $W$ is a bidiagonal pair.  By \cite[Lemma 4.2]{ItoTerWang06} and since $V$ is segregated, there exists a basis $\lbrace u_i \rbrace _{i=0}^d$ for $W$ such that 
\begin{eqnarray}
\label{uqsl2bidiag1}
x. u_i &=& q^{d-2i} u_i \qquad (0 \leq i \leq d), \\
\label{uqsl2bidiag2}
(y - q^{2i-d}I). u_i &=& (q^{-d} - q^{2i+2-d}) u_{i+1} \qquad (0 \leq i \leq d), \qquad u_{d+1}=0, \\
\label{uqsl2bidiag3}
(z - q^{2i-d}I). u_i &=& (q^{d} - q^{2i-2-d}) u_{i-1} \qquad (0 \leq i \leq d), \qquad u_{-1}=0.
\end{eqnarray}
Let $\Omega : W \rightarrow W$ denote the invertible linear operator from Lemma \ref{thm:omega}.  By Lemma \ref{thm:omega}, we have $y \, \Omega ^2 = \Omega ^2 \, x$.  Using this and (\ref{uqsl2bidiag1}), we find $y. \, \Omega ^2 u_i = q^{d-2i} \, \Omega ^2 u_i$.  Observe that $\lbrace \Omega ^2 u_i \rbrace _{i=0}^d$ is a basis for $W$, since $\Omega ^2$ is invertible.  Combining the previous two sentences, we find that the action of $y$ on $W$ is diagonalizable with eigenvalues $\lbrace q^{d-2i} \rbrace _{i=0}^d$ corresponding to eigenvectors $\lbrace \Omega ^2 u_i \rbrace _{i=0}^d$.  By Lemma \ref{thm:omega}, we have $\Omega ^2 \, y = z \, \Omega ^2$.  Using this and applying $\Omega ^ 2$ to (\ref{uqsl2bidiag2}), we have
\begin{eqnarray}
\label{uqsl2bidiag4}
z. \, \Omega ^2 u_i =  q^{2i-d}  \Omega ^2 u_i + (q^{-d} - q^{2i+2-d}) \Omega ^2 u_{i+1}  \qquad (0 \leq i \leq d).
\end{eqnarray}
Using (\ref{uqsl2bidiag4}), we have $[y,z]. \Omega ^2 u_i = q^{-2i-2} (q^2 - 1)(q^{2i+2} -1) \Omega ^2 u_{i+1}$ for $0 \leq i \leq d$.  This gives 
\begin{eqnarray}
\label{uqsl2bidiag4.5}
[y,z]^{d-2i}. \, \Omega ^2 u_i = (q^2 - 1)^{d-2i} \Big( \prod_{j=i}^{d-i-1} q^{-2j-2} (q^{2j+2} - 1) \Big) \, \Omega ^2 u_{d-i} \qquad (0 \leq i \leq d/2).
\end{eqnarray}
Define $w_i := u_{d-i}$ for $0 \leq i \leq d$.  Observe that $\lbrace w_i \rbrace _{i=0}^d$ is a basis for $W$, and that (\ref{uqsl2bidiag1}) yields 
\begin{eqnarray}
\label{uqsl2bidiag5}
x. w_i &=& q^{2i-d} w_i \qquad (0 \leq i \leq d).
\end{eqnarray}
By Lemma \ref{thm:omega}, we have $z \, \Omega = \Omega \, x$.  Using this and (\ref{uqsl2bidiag5}), we find $z \, \Omega w_i = q^{2i-d} \, \Omega w_i$ for $0 \leq i \leq d$.  Observe that $\lbrace \Omega w_i \rbrace _{i=0}^d$ is a basis for $W$, since $\Omega$ is invertible.  Combining the previous two sentences, we find that the action of $z$ on $W$ is diagonalizable with eigenvalues $\lbrace q^{2i-d} \rbrace _{i=0}^d$ corresponding to eigenvectors $\lbrace \Omega w_i \rbrace _{i=0}^d$.  From (\ref{uqsl2bidiag3}), we have 
\begin{eqnarray}
\label{uqsl2bidiag6}
(z - q^{d-2i}I). \, w_i &=& (q^{d} - q^{-2i-2+d}) w_{i+1} \qquad (0 \leq i \leq d).
\end{eqnarray}
By Lemma \ref{thm:omega}, we have $\Omega \, z = y \, \Omega$.  Using this and applying $\Omega$ to (\ref{uqsl2bidiag6}), we have 
\begin{eqnarray}
\label{uqsl2bidiag7}
y. \, \Omega w_i = q^{d-2i}  \Omega w_i + (q^d - q^{-2i-2+d}) \Omega w_{i+1} \qquad (0 \leq i \leq d).
\end{eqnarray}
Using (\ref{uqsl2bidiag7}), we find $[y,z]. \Omega w_i = -q^{-2} (q^2 - 1)(q^{2i+2} -1) \Omega w_{i+1}$ for $0 \leq i \leq d$.  This gives  
\begin{eqnarray}
\label{uqsl2bidiag7.5}
[y,z]^{d-2i}. \, \Omega w_i = (-q^{-2})^{d-2i} (q^2 - 1)^{d-2i} \Big( \prod_{j=i}^{d-i-1} (q^{2j+2} - 1) \Big) \, \Omega w_{d-i} \qquad (0 \leq i \leq d/2).
\end{eqnarray}
Recall $q$ is not a root of unity.  Hence, combining (\ref{uqsl2bidiag4}), (\ref{uqsl2bidiag4.5}), (\ref{uqsl2bidiag7}), (\ref{uqsl2bidiag7.5}), we find that the action of $y, \, z$ on $W$ is a bidiagonal pair.  Let $V = \sum_{j=0}^{r} V(d_j, 1)$ denote the direct sum decomposition of $V$ into irreducible $\uq$-submodules (since $V$ is segregated, we know $\epsilon_j = 1$ for $0 \leq j \leq r$).  Without loss of generality, $d_0 \geq d_1 \geq \cdots \geq d_r$.  For $0 \leq j \leq r$, our work above has shown that the action of $y$ (resp.~$z$) on $V(d_j, 1)$ is diagonalizable with eigenvalues $\lbrace q^{d_j -2i} \rbrace _{i=0}^{d_j}$ (resp.~$\lbrace q^{2i - d_j} \rbrace _{i=0}^{d_j}$).  Since $V$ is segregated, then either $\lbrace d_j \rbrace _{j=0}^r \subseteq 2 \Z$ or $\lbrace d_j \rbrace _{j=0}^r \subseteq 2 \Z +1$.  Combining the previous two sentences, we find that the action of $y$ (resp.~$z$) on $V$ is diagonalizable with eigenvalues $\lbrace q^{d_0 -2i} \rbrace _{i=0}^{d_0}$ (resp.~$\lbrace q^{2i - d_0} \rbrace _{i=0}^{d_0}$).  For $0 \leq i \leq d_0$, let $V_i$ (resp.~$V_i^*$) denote the eigenspace of $y$ (resp.~$z$) corresponding to $q^{d_0 -2i}$ (resp.~$q^{2i - d_0}$).   Our work above has shown $\dim(V_i) = \dim(V_{d-i})$ and $\dim(V^*_i) = \dim(V^*_{d-i})$ for $0 \leq i \leq d$.  Since (\ref{uqsl2bidiag4}), (\ref{uqsl2bidiag4.5}), (\ref{uqsl2bidiag7}), (\ref{uqsl2bidiag7.5}) hold on each $V(d_j, 1)$, we have $y, \, z, \, \lbrace V_i \rbrace _{i=0}^{d_0}, \, \lbrace \Vb_i \rbrace _{i=0}^{d_0}$ satisfy (\ref{eq:b1})--(\ref{eq:b4}).  Thus, the action of $y, \, z$ on $V$ is a bidiagonal pair with eigenvalue (resp.~dual eigenvalue) sequence $\lbrace q^{d_0-2i} \rbrace _{i=0}^{d_0}$ (resp.~$\lbrace q^{2i-d_0} \rbrace _{i=0}^{d_0}$).  Comparing Theorem \ref{thm:uqeq} and (\ref{main1.1}), we find that the base of $y, \, z$ is $q^{-2}$, which is not equal to $1$.  The previous two sentences show that $y, \, z$ is reduced.
\hfill $\Box$ \\

\section{The subspaces $W_i$}

\noindent
In this section we introduce a sequence of subspaces $\lbrace W_i \rbrace _{i=0}^d$ of $V$, and develop some of their properties.  These subspaces are the main ingredient in the  proof of Theorem \ref{thm:main3}.

\begin{definition}
\label{def:W}
\rm
With reference to Assumption \ref{assume}, define
\begin{eqnarray*}
W_i := (\Vb_0 + \cdots + \Vb_i) \cap (V_0+ \cdots +V_{d-i}) \qquad (0 \leq i \leq d).
\end{eqnarray*}
For notational convenience, let $W_{-1}:=0$ and $W_{d+1}:=0$.
\end{definition}

\noindent
The goal of this section is to prove the following theorem.

\begin{theorem}
\label{thm:Wfinal}
With reference to Definition \ref{def:decomp} and Definition \ref{def:W}, the sequence $\lbrace W_i \rbrace _{i=0}^d$ is a decomposition of $V$.
\end{theorem}

\noindent
We prove Theorem \ref{thm:Wfinal} in three steps.  First, we show that the sum $\sum_{i=0}^dW_i$ is direct.  Second, we show $V=\sum_{i=0}^dW_i$.  Finally, we show that $W_i \neq 0$ for $0 \leq i \leq d$. \\

\noindent
The arguments in this section are essentially the same as the arguments from \cite[Section 5]{Funk-Neubauer07}.  For the sake of completeness and accessibility we reproduce the arguments here in full. \\

\noindent
The following definition and the next two lemmas will be useful in proving that the sum $\sum_{i=0}^dW_i$ is direct.

\begin{definition}
\label{def:Wij}
\rm
With reference to Assumption \ref{assume}, define
\begin{eqnarray*}
W(i,j) := \biggl(\sum_{h=0}^{i}\Vb_h\biggr) \, \cap \, \biggl(\sum_{h=0}^{j}V_h\biggr) \qquad (-1 \leq i, j \leq d+1),
\end{eqnarray*}
where $\sum_{h=0}^{-1} V_h = 0$ and $\sum_{h=0}^{-1} V^*_h = 0$.
\end{definition}

\noindent
With reference to Definition \ref{def:W}, observe $W(i,d-i)=W_i$ for $0 \leq i \leq d$.

\begin{lemma}
\label{thm:AAbWij}
With reference to Assumption \ref{assume} and Definition \ref{def:Wij}, the following {\rm (i)--(iv)} hold.
\begin{enumerate}
\item[\rm (i)] $(A-\theta_jI)W(i,j) \subseteq W(i+1,j-1)$ \qquad $(0 \leq i,j \leq d)$,
\item[\rm (ii)] $(\Ab-\Tb_iI)W(i,j) \subseteq W(i-1,j+1)$ \qquad $(0 \leq i,j \leq d)$,
\item[\rm (iii)] $A W_i \subseteq W_i + W_{i+1}$ \qquad $(0 \leq i \leq d)$,
\item[\rm (iv)] $A^* W_i \subseteq W_{i-1} + W_{i}$ \qquad $(0 \leq i \leq d)$.
\end{enumerate}
\end{lemma}

\noindent
{\it Proof:} (i) Recall $V_j$ is the eigenspace for $A$ corresponding to eigenvalue $\theta_j$.  Using Definition \ref{def:Wij}, we have $(A-\theta_jI)W(i,j) \subseteq \sum_{h=0}^{j-1}V_h$.  Using Definition \ref{def:Wij} and (\ref{eq:b1}), we have $(A-\theta_jI)W(i,j) \subseteq \sum_{h=0}^{i+1}\Vb_h$.  Combining the previous two sentences, we obtain the desired result. \\
(ii) Similar to (i). \\
(iii), (iv) Immediate from (i), (ii), and Definition \ref{def:W}.
\hfill $\Box$ \\

\begin{lemma}
\label{thm:wij=0}
With reference to Definition \ref{def:Wij}, the following holds.
\begin{eqnarray*}
W(i,d-1-i)=0  \qquad (0 \leq i \leq d-1).
\end{eqnarray*}
\end{lemma}

\noindent
{\it Proof:} Define $T:=\sum_{i=0}^{d-1}W(i,d-1-i)$.  To obtain the result it suffices to show that $T=0$.  By Lemma \ref{thm:AAbWij}(ii), we find that $\Ab \, T \subseteq T$.  Recall $\Ab$ is diagonalizable on $V$, and so $\Ab$ is diagonalizable on $T$.  Also, $\Vb_j \cap T$ are the eigenspaces of $\Ab |_{T}$ for $0 \leq j \leq d$.  Thus, $T=\sum_{j=0}^{d}(\Vb_j \cap T)$ (direct sum).  Suppose, towards a contradiction, that $T \neq 0$.  Then there exists $j$ $(0 \leq j \leq d)$ such that $\Vb_j \cap T \neq 0$.  Define $t:=\min\{\,i \, |\, 0 \leq i \leq d\, , \, \Vb_i \cap T \neq 0 \,\}$ and $r:=\max\{\,i \, |\, 0 \leq i \leq d\, , \, \Vb_i \cap T \neq 0 \,\}$.  By construction $t \leq r$.  We now show that
\begin{eqnarray}
\label{numb1}
r+t \geq d.
\end{eqnarray}
If $d/2  <  t$ then (\ref{numb1}) holds, since $t \leq r$.  So now assume that $0 \leq t \leq d/2$.  Let $x \in \Vb_t \cap T$ such that $x \neq 0$.  Lemma \ref{thm:raise} gives $[A,A^*]^{d-2t}x \in \Vb_{d-t}$.  Also, by (\ref{eq:b3}), $[A,A^*]^{d-2t}x \neq 0$.  By Lemma \ref{thm:AAbWij}(i), $AT \subseteq T$.  Using this and Lemma \ref{brackpoly}(iii), $[A,A^*]^{d-2t}x \in T$.  Combining the previous four sentences, we find $\Vb_{d-t} \cap T \neq 0$.  So $d-t \leq r$, and (\ref{numb1}) holds.  Since $T \neq 0$, there exists $j$ ($0 \leq j \leq d-1$) such that $W(j, d-1-j) \neq 0$.  Define $y:=\max\{\, i \, |\, 0 \leq i \leq d-1\, , \, W(i,d-1-i) \neq 0 \,\}$.  By the definition of $T$, $T \subseteq \Vb_0+ \cdots +\Vb_y$, and so
\begin{eqnarray}
\label{numb2}
y \geq r.
\end{eqnarray}
We now show that
\begin{eqnarray}
\label{numb3}
d-y \geq t+1.
\end{eqnarray}
By definition of $y$, $W(y,d-1-y) \neq 0$.  By Definition \ref{def:Wij}, $W(y,d-1-y) \subseteq V_0+ \cdots +V_{d-1-y}$.  Combining the previous two sentences with the fact that $\sum_{i=0}^{d} V_i$ is a direct sum gives $W(y,d-1-y) \nsubseteq V_{d-y}+ \cdots +V_d$.  Therefore, $T \nsubseteq V_{d-y}+ \cdots +V_d$.  Using this and Lemma \ref{thm:Vsum}, $T \nsubseteq \Vb_{d-y}+ \cdots + \Vb_d$.  So $t < d-y$, and (\ref{numb3}) follows.  Adding (\ref{numb1}), (\ref{numb2}), and (\ref{numb3}), we find $0 \geq 1$, for a contradiction.  Thus, $T=0$, and the result follows.
\hfill $\Box$ \\

\begin{lemma}
\label{thm:directsum}
With reference to Definition \ref{def:W}, the sum $\sum_{i=0}^{d}W_i$ is direct.
\end{lemma}

\noindent
{\it Proof:}  To obtain the result it suffices to show that $(W_0+ \cdots +W_{i-1}) \cap W_i=0$  for $1 \leq i \leq d$.  By Definition \ref{def:W}, $W_0+ \cdots +W_{i-1} \subseteq \Vb_0+ \cdots + \Vb_{i-1}$ and $W_i \subseteq V_0+ \cdots +V_{d-i}$.  From this and Definition \ref{def:Wij}, we have $(W_0+ \cdots +W_{i-1}) \cap W_i \subseteq W(i-1,d-i)$.  But $W(i-1,d-i)=0$ by Lemma \ref{thm:wij=0}, and so $(W_0+ \cdots +W_{i-1}) \cap W_i=0$.
\hfill $\Box$ \\

\noindent
The following definition and the next three lemmas will be used in proving that $V=\sum_{i=0}^{d}W_i$.  \\

\noindent
Recall that $\mbox{End}(V)$ is the $\K$-algebra consisting of all linear transformations from $V$ to $V$.

\begin{definition}
\label{def:scriptD}
\rm
With reference to Assumption \ref{assume}, let $\mathcal{D}$ denote the $\K$-subalgebra of $\mbox{End}(V)$ generated by $[A,A^*]$.
\end{definition}

\noindent
We will be concerned with the following subspace of $V$.  With reference to Definition \ref{def:H_i} and Definition \ref{def:scriptD}, define $\mathcal{D}H_i :=\operatorname{span}\{\, Xh\, |\,X\in \mathcal{D}, \,h\in H_i\,\}$ for $0 \leq i \leq d/2$.

\begin{lemma}
\label{thm:scriptDAH_i}
With reference to Assumption \ref{assume} and Definition \ref{def:H_i}, the following holds.
\begin{eqnarray}
\mathcal{D}H_i=\sum_{j=0}^{d-2i} [A, A^*]^j H_i \qquad \hbox{(direct sum)} \qquad (0 \leq i \leq d/2).
\label{scriptDR^iH_i1}
\end{eqnarray}
\end{lemma}

\noindent
{\it Proof:} Define $\Delta :=\sum_{j=0}^{d-2i} [A,A^*]^j H_i$.  We first show $\mathcal{D}H_i=\Delta$.  By construction $\Delta \subseteq \mathcal{D}H_i$.  We now show that $\mathcal{D}H_i \subseteq \Delta$.  Since $\mathcal{D}$ is generated by $[A,A^*]$, and since $H_i \subseteq \Delta$, it suffices to show that $\Delta$ is $[A,A^*]$-invariant.  By (\ref{H_i1}), $[A,A^*] \Delta = \sum_{j=0}^{d-2i-1} [A,A^*]^{j+1} H_i \subseteq \Delta$.  So $\Delta$ is $[A,A^*]$-invariant, and it follows that $\mathcal{D}H_i \subseteq \Delta$.  We have now shown $\mathcal{D}H_i=\Delta$.  It remains to show that the sum $\sum_{j=0}^{d-2i} [A,A^*]^j H_i$ is direct.  By Lemma \ref{thm:raise}, we have $[A,A^*]^j H_i \subseteq V_{i+j}$ for $0 \leq j \leq d-2i$.  From this and since $\sum_{i=0}^d V_i$ is a direct sum, we find that the sum $\sum_{j=0}^{d-2i} [A,A^*]^j H_i$ is direct.  
\hfill $\Box$ \\

\begin{lemma}
\label{thm:VscriptDH_i}
With reference to Assumption \ref{assume} and Definition \ref{def:H_i},
\begin{eqnarray*}
V=\sum_{i=0}^{d/2}\mathcal{D}H_i \qquad \hbox{(direct sum)}.
\end{eqnarray*}
\end{lemma}

\noindent
{\it Proof:} By (\ref{refineV1}) and since $\lbrace V_i \rbrace _{i=0}^d$ is a decomposition of $V$, we have
\begin{eqnarray*}
V=\sum_{i=0}^{d}\sum_{j=0}^{\min(i,d-i)} [A,A^*]^{i-j} H_j \qquad \hbox{(direct sum)}.
\end{eqnarray*}
In this sum we interchange the order of summation to get
\begin{eqnarray*}
V=\sum_{i=0}^{d/2} \sum_{j=0}^{d-2i} [A,A^*]^j H_i \qquad \hbox{(direct sum)}.
\end{eqnarray*}
The result now follows by Lemma \ref{thm:scriptDAH_i}.
\hfill $\Box$ \\

\begin{lemma}
\label{thm:HW}
With reference to Definition \ref{def:H_i} and Definition \ref{def:W}, the following holds.
\begin{eqnarray*}
H_i \subseteq W_{d-i} \qquad (0 \leq i \leq d/2).
\end{eqnarray*}
\end{lemma}

\noindent
{\it Proof:}  Let $h \in H_i$, and observe by (\ref{H_i1}) that
\begin{eqnarray}
\label{HW1}
h \in V_0 + \cdots + V_i.
\end{eqnarray}
Also by (\ref{H_i1}), $[A,A^*]^{d-2i+1} h = 0$.  Using this and Lemma \ref{brackpoly}(iv), $(A^* - \theta^*_{d-i} I) \cdots (A^* - \theta^*_{i} I) h =0$.  From this we have
\begin{eqnarray}
\label{HW2}
h \in \Vb_i + \cdots + \Vb_{d-i} \subseteq \Vb_0 + \cdots + \Vb_{d-i}.
\end{eqnarray}
Combining (\ref{HW1}) and (\ref{HW2}) with Definition \ref{def:W} gives $h \in W_{d-i}$, and the result follows.
\hfill $\Box$ \\

\begin{lemma}
\label{thm:wsumv}
With reference to Definition \ref{def:W},
\begin{eqnarray*}
\label{wdirectsumv1}
V=\sum_{i=0}^{d}W_i.
\end{eqnarray*}
\end{lemma}

\noindent
{\it Proof:} Define $V' :=\sum_{i=0}^{d}W_i$.  We show $V=V'$.  By construction $V' \subseteq V$.  We now show that $V \subseteq V'$.  By Lemma \ref{thm:AAbWij}(iii),(iv) we have $[A,A^*] W_i \subseteq W_{i-1} + W_i + W_{i+1}$ for $0 \leq i \leq d$.  Thus, $[A,A^*] V' \subseteq V'$.  By this and Definition \ref{def:scriptD}, $\mathcal{D}V' \subseteq V'$.  By Lemma \ref{thm:HW}, $H_j \subseteq V'$ for $0 \leq j \leq d/2$.  The previous two sentences give $\mathcal{D}H_j \subseteq V'$ for $0 \leq j \leq d/2$.  From this and Lemma \ref{thm:VscriptDH_i}, $V \subseteq V'$.  We have now shown $V=V'$.
\hfill $\Box$ \\

\begin{corollary}
\label{thm:WVVbfinal}
With reference to Assumption \ref{assume} and Definition \ref{def:W}, the following \\ {\rm (i)--(iii)} hold.
\begin{enumerate}
\item[\rm (i)] $W_0+ \cdots +W_i= \Vb_0+ \cdots + \Vb_i$  \qquad $(0 \leq i \leq d)$,
\item[\rm (ii)] $\dim(W_i)= \dim(V_i ^*)$ \qquad $(0 \leq i \leq d)$,
\item[\rm (iii)] $W_i \neq 0$ \qquad $(0 \leq i \leq d)$.
\end{enumerate}
\end{corollary}

\noindent
{\it Proof:} (i) Define $\Delta := W_0+ \cdots +W_i$ and $\Gamma := \Vb_0+ \cdots + \Vb_i$.  We show $\Delta = \Gamma$.  By Definition \ref{def:W}, $\Delta \subseteq \Gamma$, and so $\dim(\Delta) \leq \dim(\Gamma)$.  Thus, to obtain the result it suffices to show $\dim(\Delta) = \dim(\Gamma)$.  Suppose, towards a contradiction, that $\dim(\Delta) < \dim(\Gamma)$.  Then by Lemma \ref{thm:directsum} and since $\sum_{i=0}^d \Vb_i$ is a direct sum, we have
\begin{eqnarray}
\label{stop1}
\sum_{h=0}^{i} \dim(W_h) < \sum_{h=0}^{i} \dim(\Vb_h).
\end{eqnarray}
By Definition \ref{def:W}, $\sum_{h=i+1}^dW_h \subseteq \sum_{h=0}^{d-i-1}V_h$.  From this, Lemma \ref{thm:samedim}, Lemma \ref{thm:directsum}, and since $\sum_{i=0}^d V_i$ is a direct sum, we have
\begin{eqnarray}
\label{stop2}
\sum_{h=i+1}^{d} \dim(W_{h}) \leq \sum_{h=i+1}^{d} \dim(V_{h}).
\end{eqnarray}
By Lemma \ref{thm:directsum} and Lemma \ref{thm:wsumv}, we have
\begin{eqnarray}
\label{stop3}
\dim(V) = \sum_{h=0}^{d} \dim(W_h).
\end{eqnarray}
Adding (\ref{stop1})--(\ref{stop3}) we find that $\dim(V)$ $<$ $\sum_{h=0}^{i} \dim(\Vb_h) + \sum_{h=i+1}^{d} \dim(V_h)$.  Using this and Lemma \ref{thm:samedim}, we find that $\dim(V) < \sum_{h=0}^{d} \dim(V_h)$.  From this and since $\lbrace V_i \rbrace _{i=0}^d$ is a decomposition of $V$, we have $\dim (V) < \dim (V)$, for a contradiction.  Thus, $\dim(\Delta) = \dim(\Gamma)$, and the result follows.  \\
(ii) By (i), Lemma \ref{thm:directsum}, and since $\sum_{i=0}^d V_i ^*$ is a direct sum, we find $\sum_{h=0}^{i} \dim(W_h)=\sum_{h=0}^{i} \dim(V_h ^*)$ for $0 \leq  i \leq d$.  The result follows immediately from this. \\
(iii) For $0 \leq i \leq d$, $V_i^* \neq 0$, since $\{ V_i ^* \}_{i=0} ^d$ is a decomposition of $V$.  From this and (ii), the result follows immediately.
\hfill $\Box$ \\

\noindent
We have now reached the goal of this section. \\

\noindent
{\it Proof of Theorem \ref{thm:Wfinal}:}  Combining Lemma \ref{thm:directsum}, Lemma \ref{thm:wsumv}, and Corollary \ref{thm:WVVbfinal} (iii) we immediately obtain Theorem \ref{thm:Wfinal}.
\hfill $\Box$ \\

\section{The proof that bidiagonal pairs act as equitable pairs}

In this section we prove Theorem \ref{thm:main3}. \\

\noindent
Throughout this section we adopt Assumption \ref{assume}, along with the additional assumption that $A, \, A^*$ is reduced.  Also, throughout this section $b$ will denote the base of $A, \, A^*$.  \\

\noindent
We now introduce a linear transformation which will be used in the proof of Theorem \ref{thm:main3}(i).

\begin{definition}
\label{def:B1}
\rm
With reference to Definition \ref{def:W}, define the following linear transformation.  Let $B :V\to V$ be the linear transformation such that for $0 \leq i \leq d$, $W_i$  is an eigenspace for $B$ with eigenvalue $2i-d$.
\end{definition}

\begin{lemma}
\label{thm:sigmatauB}
Suppose that $b=1$.  Then with reference Definition \ref{def:B1}, the following hold.
\begin{eqnarray}
\label{A1Ab1B11}
A \Ab - \Ab A &=& 2 A + 2 \Ab, \\
\label{A1Ab1B12}
\Ab B - B \Ab &=& 2 \Ab + 2 B, \\
\label{A1Ab1B13}
B A - A B &=& 2 B + 2 A.
\end{eqnarray}
\end{lemma}

\noindent
{\it Proof:}  Since $A, \, A^*$ is reduced, we have 
\begin{eqnarray}
\label{sigmatauB1}
\theta_i = 2i-d, \qquad \qquad  \theta_i ^* = d-2i \qquad (0 \leq i \leq d).
\end{eqnarray}
First asume $d=0$.  Then $\theta_0 = 0, \theta_0 ^* = 0$, and $V = V_0 = V_0 ^* = W_0$.  This gives $A=0$, $A^* = 0$, and $B=0$.  So (\ref{A1Ab1B11})--(\ref{A1Ab1B13}) hold.  \\ \\
Now assume that $d \geq 1$.  Combining (\ref{sigmatauB1}) with Theorem \ref{thm:main1} and Lemma \ref{thm:relrec}, we obtain (\ref{A1Ab1B11}).  We now show (\ref{A1Ab1B12}).  Using (\ref{sigmatauB1}) and Lemma \ref{thm:AAbWij}(ii) (with $j=d-i$), we have $(A^* - (d-2i) I) W_i \subseteq W_{i-1}$ for $0 \leq i \leq d$.  This and Definition \ref{def:B1} yields $(B - (2i-2-d) I)(A^* - (d-2i) I) W_i = 0$ for $0 \leq i \leq d$.  Hence,
\begin{eqnarray*}
(BA^* - (d-2i) B - A^* (2i-d) + 2A^* + (d-2i)(2i-d) I + 2(2i-d) I ) W_i = 0.
\end{eqnarray*}
From this and Definition \ref{def:B1}, we find $- (A^*B - BA^* -2A^* -2B) W_i =0$.  Thus, we obtain (\ref{A1Ab1B12}), since $\{ W_i \}_{i=0}^d$ is a decomposition of $V$.  The proof of (\ref{A1Ab1B13}) is similar to the proof of (\ref{A1Ab1B12}).
\hfill $\Box$ \\

\noindent
We now introduce a linear transformation which will be used in the proof of Theorem \ref{thm:main3}(ii).

\begin{definition}
\label{def:B'}
\rm
Let $q$ denote a nonzero scalar in $\K$, which is not a root of unity.  With reference to Definition \ref{def:W}, define the following linear transformation.  Let $B' :V\to V$ be the linear transformation such that for $0 \leq i \leq d$, $W_i$  is an eigenspace for $B'$ with eigenvalue $q^{d-2i}$.
\end{definition}

\begin{lemma}
\label{thm:sigma'tau'B'}
Suppose that $b \neq 1$.  Then with reference to Note \ref{bnot1} and Definition \ref{def:B'}, the following hold.
\begin{eqnarray}
\label{AqAbqBq1}
\frac{q A \Ab - q^{-1} \Ab A} {q-q^{-1}} &=& I, \\
\label{AqAbqBq2}
\frac{q \Ab B' - q^{-1} B' \Ab} {q-q^{-1}} &=& I, \\
\label{AqAbqBq3}
\frac{q B' A - q^{-1} A B'} {q-q^{-1}} &=& I.
\end{eqnarray}
\end{lemma}

\noindent
{\it Proof:}  Since $A, \, A^*$ is reduced, we have 
\begin{eqnarray}
\label{sigmatauB'}
\theta_i = q^{d-2i}, \qquad \qquad  \theta_i ^* = q^{2i-d} \qquad (0 \leq i \leq d).
\end{eqnarray}
Combining (\ref{sigmatauB'}) with Theorem \ref{thm:main1} and Lemma \ref{thm:relrec} we obtain (\ref{AqAbqBq1}).  We now show (\ref{AqAbqBq2}).  Using (\ref{sigmatauB'}) and Lemma \ref{thm:AAbWij}(ii) (with $j=d-i$), we have $(A^* - q^{2i-d} I) W_i \subseteq W_{i-1}$ for $0 \leq i \leq d$.  This and Definition \ref{def:B'} yields $(B' - q^{d-2i+2} I)(A^* - q^{2i-d} I) W_i = 0$ for $0 \leq i \leq d$.  Hence,
\begin{eqnarray*}
(B'A^* - q^{2i-d} B' - q^2 A^* q^{d-2i} + q^2 I ) W_i = 0.
\end{eqnarray*}
From this and Definition \ref{def:B'}, we find $-q (qA^*B' - q^{-1} B'A^* - (q-q^{-1}) I) W_i =0$.  Thus, we obtain (\ref{AqAbqBq2}), since $\{ W_i \}_{i=0}^d$ is a decomposition of $V$.  The proof of (\ref{AqAbqBq3}) is similar to the proof of (\ref{AqAbqBq2}).
\hfill $\Box$ \\

\noindent
The following lemma will be used in the proof of Theorem \ref{thm:main3}(i).

\begin{lemma}
\label{thm:sl2unique}
Let $h, \, e, \, f$ be as in Definition \ref{def:usl2}, and let $X, \, Y, \, Z$ be an equitable basis for $\SL$.  Let $V$ denote a vector space over $\K$ with finite positive dimension.  Suppose that there are two $\SL$-module structures on $V$.  Then the following holds. \\ \\
Assume that the actions of $h$ (resp.~$Z$) on $V$ given by the two module structures agree.  Assume that the actions of $f$ (resp.~$Y$) on $V$ given by the two module structures agree.  Then the actions of $e$ (resp.~$X$) on $V$ given by the two module structures agree.
\end{lemma}

\noindent
{\it Proof:} First we prove the result involving $h, e, f$.  Let $E_1:V \to V$ (resp.~$E_2:V \to V$) denote the action of $e$ on $V$ given by the first (resp.~second) module structure.  We show $(E_1-E_2)V=0$.  Using Lemma \ref{thm:compred} and referring to the first module structure, $V$ is the direct sum of irreducible $\SL$-submodules.  Let $W$ be one of the irreducible submodules in this sum.  It suffices to show $(E_1-E_2)W=0$.  By Lemma \ref{thm:sl2mods}, there exists a nonnegative integer $d$  such that $W$ is isomorphic to $V(d)$.  Therefore, the eigenvalues for $h$ on $W$ are $d-2i$ $(0 \leq i \leq d)$, and $\dim(W)=d+1$.  Let $w \in W$ be an eigenvector for $h$ with eigenvalue $d$.  By Lemma \ref{thm:sl2mods}, $\lbrace f^{i}. w \rbrace _{i=0}^d$ is a basis for $W$.  We show by induction that $(E_1-E_2)f^{i}. w=0$ for $0 \leq i \leq d$.  First assume that $i=0$.  By Lemma \ref{thm:sl2mods}, $E_1 w=0$.  Also by Lemma \ref{thm:sl2mods}, $f^{d+1}. w=0$, and so $E_2 w=0$.  We have now shown $(E_1 - E_2)w=0$.  Next assume that $i \geq 1$.  By induction we have
\begin{eqnarray}
\label{sl2unique1}
(E_1 - E_2)f^{i-1}. w=0.
\end{eqnarray}
By Definition \ref{def:usl2} and since the actions of $h$ (resp.~$f$) on $V$ given by the two module structures agree, we have $[E_1 , f] = [E_2 , f]$. Using this we have
\begin{eqnarray}
\label{sl2unique2}
(E_1 f - E_2 f) f^{i-1}. w = f. (E_1 - E_2) f^{i-1}. w
\end{eqnarray}
Combining (\ref{sl2unique1}) and (\ref{sl2unique2}) gives $(E_1 - E_2)f^{i}. w=0$.  We have now shown $(E_1 - E_2)W=0$, and so $(E_1 - E_2)V=0$.  Thus, the actions of $e$ on $V$ given by the two module structures agree. \\ \\
We now prove the result involving $X, \, Y, \, Z$.  Identify the copy of $\SL$ given in Definition \ref{def:usl2} with the copy given in Theorem \ref{thm:ueq}, via the automorphism given in Theorem \ref{thm:ueq}.  Under this automorphism, $Z$ is mapped to $h$, $Y$ is mapped to $-2f-h$, and $X$ is mapped to $2e-h$.  From this and the result involving $h, \, e, \, f$, we find that the actions of $X$ on $V$ given by the two module structures agree.
\hfill $\Box$ \\

\noindent
The following lemma will be used in the proof of Theorem \ref{thm:main3}(ii).

\begin{lemma}
\label{thm:uqsl2unique}
Let $k, \, e, \, f$ be as in Defnition \ref{def:uqsl2}, and let $x^{\pm 1}, \, y, \, z$ be the equitable generators for $\uq$.  Let $V$ denote a vector space over $\K$ with finite positive dimension.  Suppose that there are two $\uq$-module structures on $V$.  Then the following holds. \\ \\
Assume that the actions of $k$ (resp.~$y$) on $V$ given by the two module structures agree.  Assume that the actions of $f$ (resp.~$z$) on $V$ given by the two module structures agree.  Then the actions of $e$ (resp.~$x^{\pm 1}$) on $V$ given by the two module structures agree.
\end{lemma}

\noindent
{\it Proof:} The result involving $k, \, e, \, f$ is proven in \cite[Lemma 9.8]{Funk-Neubauer07}. \\ \\
We now prove the result involving $x^{\pm 1}, \, y, \, z$.  By \cite[Corollary 4.5]{ItoTerWang06}, the action of $y$ on $V$ is invertible.  Let $y^{-1}$ denote the inverse of the action of $y$ on $V$.  Identify the $\uq$-module structure on $V$ given by $k^{\pm1}, \, e, \, f$ with the $\uq$-module structure on $V$ given by $x, \, y^{\pm1}, \, z$, via the following isomorphism:
\begin{eqnarray*}
k^{\pm1} &\rightarrow& y^{\pm1}, \\
f &\rightarrow& (z - y^{-1})(q-q^{-1}), \\
e &\rightarrow& (1 - yx)q^{-1}(q-q^{-1}).
\end{eqnarray*}
From this and the result involving $k, \, e, \, f$, we have that the actions of $x^{\pm 1}$ on $V$ given by the two module structures agree.
\hfill $\Box$ \\

\noindent
We are now ready to prove Theorem \ref{thm:main3}. \\

\noindent
{\it Proof of Theorem \ref{thm:main3}:}  \\
(i) Suppose that $b = 1$.  Let $Y, \, Z$ denote an equitable pair in $\SL$.  Let $B$ be as in Definition \ref{def:B1}.  Comparing Theorem \ref{thm:ueq} and Lemma \ref{thm:sigmatauB}, we find there exists a $\SL$-module structure on $V$ such that $(X - B)V=0$, $(Y - A)V=0$, and $(Z - \Ab)V=0$.  By Lemma \ref{thm:sl2unique}, the action of $X$ on $V$ is uniquely determined by the actions of $Y$ and $Z$ on $V$.  Thus, the action of $X$ on $V$ is uniquely determined by $A, \, A^*$, since $(Y - A)V=0$ and $(Z - \Ab)V=0$.  We now check that this $\SL$-module structure on $V$ is segregated.  The dual eigenvalue sequence of $A, \, A^*$ is $\lbrace d-2i \rbrace _{i=0}^d$, since $A, \, A^*$ is reduced.  Using the $\SL$ automorphism from Theorem \ref{thm:ueq}, we find $(h - Z)V = 0$.  Combining the previous two sentences with $(Z - \Ab)V=0$, we find that the action of $h$ on $V$ has eigenvalues $\lbrace d-2i \rbrace _{i=0}^d$.  So if $d$ is even (resp.~odd) then $V=V_{\hbox{even}}$ (resp.~$V=V_{\hbox{odd}}$).  Thus, the $\SL$-module structure on $V$ is segregated. \\ \\
(ii) Suppose that $b \neq 1$.  Let $y, \, z$ denote an equitable pair in $\uq$.  Let $B'$ be as in Definition \ref{def:B'}.  Comparing Theorem \ref{thm:uqeq} and Lemma \ref{thm:sigma'tau'B'}, we find there exists a $\uq$-module structure on $V$ such that $(x - B')V=0$, $(x^{-1} - B'^{-1})V=0$, $(y - A)V=0$, and $(z - \Ab)V=0$. By Lemma \ref{thm:uqsl2unique}, the actions of $x^{\pm1}$ on $V$ are uniquely determined by the actions of $y$ and $z$ on $V$.  Thus, the actions of $x^{\pm1}$ on $V$ are uniquely determined by $A, \, A^*$, since $(y - A)V=0$ and $(z - \Ab)V=0$.  We now check this $\uq$-module structure on $V$ is segregated.  Using the $\uq$ automorphism from Theorem \ref{thm:uqeq}, we have
\begin{eqnarray}
\label{kx}
(k - x)V = 0.
\end{eqnarray}
Let $\Omega : V \rightarrow V$ be the invertible linear operator from Lemma \ref{thm:omega}, and recall
\begin{eqnarray}
\label{omeg}
(\Omega ^{-1} z \Omega - x)V = 0.
\end{eqnarray} 
The dual eigenvalue sequence of $A, \, A^*$ is $\lbrace q^{2i-d} \rbrace _{i=0}^d$, since $A, \, A^*$ is reduced.  Combining this with (\ref{kx}), (\ref{omeg}), and $(z - \Ab)V=0$, we have $k. \Omega ^{-1} v = q^{2i-d} \Omega ^{-1} v$ for $0 \neq v \in V_i ^*$.  Thus, the action of $k$ on $V$ has eigenvalues $\lbrace q^{2i-d} \rbrace _{i=0}^d$.  So if $d$ is even (resp.~odd) then $V=V_{\hbox{even}}^1$ (resp.~$V=V_{\hbox{odd}}^1$).  Thus, the $\uq$-module structure on $V$ is segregated.
\hfill $\Box$ \\

\section{The proof of the classification theorem}

\noindent
In this section we prove Theorem \ref{thm:class}.  Lemma \ref{thm:isompa} will be used in the proof of Theorem \ref{thm:class}.  So we now prove Lemma \ref{thm:isompa}, namely that two bidiagonal pairs are isomorphic exactly when their parameter arrays are equal. \\

\noindent
{\it Proof of Lemma \ref{thm:isompa}:}  \\
Let $A, \, A^*$ and $B, \, B^*$ denote bidiagonal pairs over $\K$.  Let $V$ (resp.~$\widetilde{V}$) denote the vector space underlying $A, \, A^*$ (resp.~$B, \, B^*$).  Let
\begin{eqnarray}
\label{isompa1}
(\{ \theta_i \}_{i=0}^d; \{  \theta^*_i \}_{i=0}^d; \{ \rho_i \}_{i=0}^d)
\end{eqnarray}
denote the parameter array of $A, \, A^*$.  For $0 \leq i \leq d$, let $V_i$ denote the eigenspace of $A$ corresponding to $\theta_i$.  Observe $\{  V_i \}_{i=0}^d$ is the standard ordering of the eigenspaces of $A$.  \\

\noindent
$(\Longrightarrow)$:   Suppose that $A, \, A^*$ and $B, \, B^*$ are isomorphic, and let $\mu : V \rightarrow \widetilde{V}$ denote an isomorphism of bidiagonal pairs from $A, \, A^*$ to $B, \, B^*$.  We show that (\ref{isompa1}) is the parameter array of $B, \, B^*$.  We first show that $\theta_i$ is an eigenvalue of $B$ for $0 \leq i \leq d$.  Let $0 \neq v \in V_i$.  Observe $\mu v \neq 0$, since $\mu$ is injective.  Definition \ref{def:isom} gives $B \mu v = \mu A v = \theta_i \mu v$, and so $\theta_i$ is an eigenvalue of $B$.  For $0 \leq i \leq d$, let $\widetilde{V}_i$ denote the eigenspace of $B$ corresponding to $\theta_i$.  We now show 
\begin{eqnarray}
\label{isompa2}
\widetilde{V}_i = \mu (V_i) \qquad (0 \leq i \leq d).
\end{eqnarray}
First we show $\widetilde{V}_i \subseteq \mu (V_i)$.  Let $x \in \widetilde{V}_i$.  Since $\mu$ is surjective, there exists $v \in V$ such that $\mu v = x$.  Definition \ref{def:isom} gives $\mu A v = B \mu v = \theta_i x = \mu \theta_i v$.  So $A v = \theta_i v$, since $\mu$ is injective.  Thus, $v \in V_i$, and so $x \in \mu (V_i)$.  We have now shown that $\widetilde{V}_i \subseteq \mu (V_i)$.  Next we show that $\mu (V_i) \subseteq \widetilde{V}_i$.  Let $x \in \mu (V_i)$, and let $v \in V_i$ with $x = \mu v$.  Definition \ref{def:isom} gives $B x = \mu A v = \mu \theta_i v = \theta_i x$, and so $x \in \widetilde{V}_i$.  We have now shown that $\mu (V_i) \subseteq \widetilde{V}_i$, and so (\ref{isompa2}) holds.  We now show that $B^* \widetilde{V}_i \subseteq \widetilde{V}_i + \widetilde{V}_{i+1}$ for $0 \leq i \leq d$.  Let $x \in \widetilde{V}_i$.  By Definition \ref{def:isom} and (\ref{isompa2}), there exists $v \in V_i$ such that $B^* x = \mu A^* v$.  From this, (\ref{eq:b2}), and (\ref{isompa2}), we find that $B^* x \in \mu (V_i + V_{i+1}) = \widetilde{V}_i + \widetilde{V}_{i+1}$.  We have now shown that $B^* \widetilde{V}_i \subseteq \widetilde{V}_i + \widetilde{V}_{i+1}$ for $0 \leq i \leq d$.  So $\{ \widetilde{V}_i \}_{i=0}^d$  is the standard ordering of the eigenspaces of $B$.  Therefore,  $\{ \theta_i \}_{i=0}^d$ is the eigenvalue sequence of $B, \, B^*$.  From this and Note \ref{A^*A}, we have $\{ \theta^*_i \}_{i=0}^d$ is the dual eigenvalue sequence of $B, \, B^*$.  We now show that $\{ \rho_i \}_{i=0}^d$ is the shape of $B, \, B^*$.  From (\ref{isompa2}), and since $\mu$ is a vector space isomorphism, we find that $\mu |_{V_i} : V_i \rightarrow \widetilde{V}_i$ is a vector space isomorphism for $0 \leq i \leq d$.  So $\rho_i = \dim(V_i) = \dim(\widetilde{V}_i) $ for $0 \leq i \leq d$.  Thus, $\{ \rho_i \}_{i=0}^d$ is the shape of $B, \, B^*$.  We have now shown that (\ref{isompa1}) is the parameter array of $B, \, B^*$, and so the parameter array of $A, \, A^*$ equals the parameter array of $B, \, B^*$.\\

\noindent
$(\Longleftarrow)$:  Suppose that the parameter array of $A, \, A^*$ equals the parameter array of $B, \, B^*$, so that (\ref{isompa1}) is also the parameter array of $B, \, B^*$.  We show that $A, \, A^*$ and $B, \, B^*$ are isomorphic.  We break the argument into the following four steps.  In step 1 we construct a certain useful basis $\{ v_{i,j,k} \}$ for $V$.  In step 2 we construct a certain useful basis $\{ \widetilde{v}_{i,j,k} \}$ for $\widetilde{V}$.  In step 3 we use the bases from step 1 and step 2 to construct a vector space isomorphism $\mu$ between $V$ and $\widetilde{V}$.  In step 4 we show that $\mu$ is an isomorphism of bidiagonal pairs from $A, \, A^*$ to $B, \, B^*$.  \\ \\
Step 1:  For $0 \leq i \leq d/2$, let $H_i$ denote the subspace of $V$ defined in (\ref{H_i1}), and let $h_i := \dim(H_i)$.  If $h_i \neq 0$, let $\{ v_{i,j} \}_{j=1}^{h_i}$ denote a basis for $H_i$.  If $h_i = 0$, we interpret $\{ v_{i,j} \}_{j=1}^{h_i}$ as the empty set.  For the remainder of this proof, the indices $i, j, k$ will always satisfy $0 \leq i \leq d/2$, $1 \leq j \leq h_i$, and $0 \leq k \leq d-2i$.  Define $v_{i,j,k} := [A,A^*] ^k v_{i, j}$, and $v_{i,j,d-2i+1} := 0$.  Observe $v_{i,j,k} \in [A,A^*] ^k H_i \subseteq V_{i+k}$.  By Lemma \ref{thm:injsurj}, $[A, A^*] ^k |_{V_i} : V_i \rightarrow V_{i+k}$ is an injection.  So $[A, A^*] ^k |_{H_i} : H_i \rightarrow [A,A^*] ^k H_i$ is a bijection, and so $\bigcup_j \{ v_{i,j,k} \}$ is a basis for $[A, A^*] ^k H_i$.  This and Lemma \ref{thm:scriptDAH_i} give that $\bigcup_{j,k} \{ v_{i,j,k} \}$ is a basis for $\mathcal{D} H_i$.  This and Lemma \ref{thm:VscriptDH_i} give that $\bigcup_{i,j,k} \{ v_{i,j,k} \}$ is a basis for $V$.  \\ \\
Step 2:  Let $\widetilde{V_i}$ denote the eigenspace of $B$ corresponding to $\theta_i$.  \\ Let $\widetilde{H}_i :=\{\,v \in \widetilde{V}_i\, |\, [B, B^*]^{d-2i+1}v=0\,\}$, and let $\widetilde{h}_i := \dim(\widetilde{H}_i)$.   Before we construct the basis $\{ \widetilde{v}_{i,j,k} \}$, we first show that 
\begin{eqnarray}
\label{hsequal}
h_i = \widetilde{h}_i.
\end{eqnarray}
By (\ref{refineVbi2}) we find $\dim(V_i) = \dim([A,A^*] V_{i-1}) + \dim(H_i)$.  From this and Lemma \ref{thm:injsurj}, we find $\dim(V_i) - \dim(V_{i-1}) = \dim(H_i)$.  Similarly, applying (\ref{refineVbi2}) and Lemma \ref{thm:injsurj} to $B, \, B^*$, we find $\dim(\widetilde{V}_i) - \dim(\widetilde{V}_{i-1}) = \dim(\widetilde{H}_i)$.  The parameter array of $A, \, A^*$ equals the parameter array of $B, \, B^*$, and so $\dim(V_i) = \dim(\widetilde{V}_i)$.  Combining the previous three sentences, we obtain (\ref{hsequal}).  If $h_i = \widetilde{h}_i \neq 0$, let $\{ \widetilde{v}_{i,j} \}_{j=1}^{h_i}$ denote a basis for $\widetilde{H}_i$.  Define $\widetilde{v}_{i,j,k} := [B,B^*] ^k \widetilde{v}_{i, j}$, and $\widetilde{v}_{i,j,d-2i+1} := 0$.  Observe $\widetilde{v}_{i,j,k} \in [B,B^*] ^k \widetilde{H}_i \subseteq \widetilde{V}_{i+k}$.  As in step 1, using Lemma \ref{thm:injsurj}, Lemma \ref{thm:scriptDAH_i}, and Lemma \ref{thm:VscriptDH_i}, we find that $\bigcup_{i,j,k} \{ \widetilde{v}_{i,j,k} \}$ is a basis for $\widetilde{V}$.  \\ \\
Step 3:  Define the linear transformation $\mu : V \rightarrow \widetilde{V}$ as follows.  Let $\mu (v_{i,j,k}) := \widetilde{v}_{i,j,k}$, and extend $\mu$ linearly to $V$.  From (\ref{hsequal}) we see that $\mu$ is well defined.  Since $\bigcup_{i,j,k} \{\widetilde{v}_{i,j,k} \}$ is a basis for $\widetilde{V}$, we find that $\mu$ is a vector space isomorphism.  \\ \\
Step 4:  We now show $\mu A = B \mu$.  Let $v \in V$.  Since $\bigcup_{i,j,k} \{ v_{i,j,k} \}$ is a basis for $V$, there exist $c_{i,j,k} \in \K$ such that $v = \sum_{i,j,k} \, c_{i,j,k} \, v_{i,j,k}$.  For notational convenience, let $c_{i,j,-1} := 0$.  Since $v_{i,j,k} \in V_{i+k}$ and $\widetilde{v}_{i,j,k} \in \widetilde{V}_{i+k}$, we have
\begin{eqnarray*}
\mu \, A \, v &=& \mu \, \big( \sum_{i,j,k} \, c_{i,j,k} \, \theta_{i+k} \, v_{i,j,k} \big)  \\
&=& \sum_{i,j,k} \, c_{i,j,k} \, \theta_{i+k} \, \widetilde{v}_{i,j,k} \\
&=&  \sum_{i,j,k} \, c_{i,j,k} \, B \, \widetilde{v}_{i,j,k} \\
&=& B \, \mu \, v.
\end{eqnarray*}
Thus, $\mu A = B \mu$.  We now show $\mu A^* = B^* \mu$.   In order to do this, we first show
\begin{eqnarray}
\label{isompa3}
A^* \, v_{i,j,k} = \theta^*_{i+k} \, v_{i,j,k} + (\theta_{i+k+1} - \theta_{i+k})^{-1} \,  v_{i,j,k+1}, \\
\label{isompa4}
B^* \, \widetilde{v}_{i,j,k} = \theta^*_{i+k} \, \widetilde{v}_{i,j,k} + (\theta_{i+k+1} - \theta_{i+k})^{-1} \, \widetilde{v}_{i,j,k+1}.
\end{eqnarray}
Lemma \ref{brackpoly}(iv) gives $v_{i,j,k} = \prod_{s=i}^{k+i-1} (\theta_{s+1} - \theta_s)(A^* - \theta^*_s I) v_{i,j}$.  This gives $(\theta_{i+k+1} - \theta_{i+k})(A^* - \theta^*_{i+k} I) v_{i,j,k} = v_{i,j,k+1}$, and we obtain (\ref{isompa3}). ÊApplying Lemma \ref{brackpoly}(iv) to $B, \, B^*$, and since (\ref{isompa1}) is the parameter array of $B, \, B^*$, we find $\widetilde{v}_{i,j,k} = \prod_{s=i}^{k+i-1} (\theta_{s+1} - \theta_s)(B^* - \theta^*_s I) \widetilde{v}_{i,j}$.  This gives $(\theta_{i+k+1} - \theta_{i+k})(B^* - \theta^*_{i+k} I) \widetilde{v}_{i,j,k} = \widetilde{v}_{i,j,k+1}$, and we obtain (\ref{isompa4}).  Using (\ref{isompa3}) and (\ref{isompa4}), we have
\begin{eqnarray*}
\mu \, A^* \, v &=& \mu \, \big( \sum_{i,j,k} \, c_{i,j,k} \, \big( \theta^*_{i+k} \, v_{i,j,k} + (\theta_{i+k+1} - \theta_{i+k})^{-1} \,  v_{i,j,k+1} \big) \, \big)  \\
&=& \mu \, \big( \sum_{i,j,k} \, \big( c_{i,j,k} \, \theta^*_{i+k} + c_{i,j,k-1} \, (\theta_{i+k} - \theta_{i+k-1})^{-1} \big) \, v_{i,j,k} \, \big) \\
&=& \sum_{i,j,k} \, \big( c_{i,j,k} \, \theta^*_{i+k} + c_{i,j,k-1} \, (\theta_{i+k} - \theta_{i+k-1})^{-1} \big) \, \widetilde{v}_{i,j,k} \\
&=& \sum_{i,j,k} \, c_{i,j,k} \, \big( \theta^*_{i+k} \, \widetilde{v}_{i,j,k} + (\theta_{i+k+1} - \theta_{i+k})^{-1} \,  \widetilde{v}_{i,j,k+1} \big)\\
&=& \sum_{i,j,k} \, c_{i,j,k} \, B^* \, \widetilde{v}_{i,j,k} \\
&=& B^* \, \mu \, v.
\end{eqnarray*}
Thus, $\mu A^* = B^* \mu$.  We have now shown that $\mu$ is an isomorphism of bidiagonal pairs from $A, \, A^*$ to $B, \, B^*$.  So $A, \, A^*$ and $B, \, B^*$ are isomorphic.
\hfill $\Box$ \\

 \noindent
The following lemma will be used in the proof of Theorem \ref{thm:class}.

\begin{lemma}
\label{thm:bijcorrsl2}
Let $d$ denote a nonnegative integer, and let $\{ \rho_i \}_{i=0}^d$ denote a sequence of scalars satisfying Theorem \ref{thm:class}(iii)--(v).  Then there exists a unique (up to isomorphism) segregated $\SL$-module (resp.~$\uq$-module) $V$ satisfying the following {\rm(i)},{\rm(ii)}.
\begin{enumerate}
\item[\rm (i)] The action of $h$ (resp.~$k$) on $V$ is diagonalizable with eigenvalues $\{ d-2i \}_{i=0}^d$ (resp.~$\{ q^{d-2i} \}_{i=0}^d$).
\item[\rm (ii)] For $0 \leq i \leq d$, $\rho_i  = \dim(U_i)$, where $U_i$ is the eigenspace for the action of $h$ (resp.~$k$) on $V$ corresponding to the eigenvalue $d-2i$ (resp.~$q^{d-2i}$).
\end{enumerate}
\end{lemma}

\noindent
{\it Proof:}   First we prove the result involving $\SL$.  For notational convenience, let $\rho_{-1} := 0$.  Define $d_j := d-2j$ and $m_j := \rho_j - \rho_{j-1}$ for $0 \leq j \leq d/2$.  Let $V(d_j)$ denote the finite-dimensional $\SL$-module from Lemma \ref{thm:sl2mods}.  Define $V := \sum_{j=0}^{d/2} V(d_j)^{+m_j}$, where $V(d_j)^{+m_j}$ denotes $V(d_j) + \cdots + V(d_j)$ ($m_j$ times).  By construction, either $\{ d_j \}_{j=0}^{d/2} \subseteq 2 \Z$ or $\{ d_j \}_{j=0}^{d/2} \subseteq 2 \Z + 1$.  From this, Lemma \ref{thm:sl2mods}, and Definition \ref{def:pure}, we find that the $\SL$-module $V$ is segregated.  By construction, $d=d_0 > d_1 > \cdots > d_{d/2}$.  From this and Lemma \ref{thm:sl2mods}, we find that the action of $h$ on $V$ is diagonalizable with eigenvalues $\{ d-2i \}_{i=0}^d$.  For $0 \leq i \leq d/2$, $\dim(U_i) = m_0 + m_1 + \cdots + m_i = \rho_i$.  For $d/2 < i \leq d$, $\dim(U_i) = m_0 + m_1 + \cdots + m_{d-i} = \rho_{d-i} = \rho_i$.  Therefore, $\dim(U_i) = \rho_i$ for $0 \leq i \leq d$.  By construction, $d$ and the sequence $\{ \rho_i \}_{i=0}^d$ uniquely determine the sequences  $\{ d_j \}_{j=0}^{d/2}$, $\{ m_j \}_{j=0}^{d/2}$.  Also, the sequences  $\{ d_j \}_{j=0}^{d/2}$, $\{ m_j \}_{j=0}^{d/2}$ uniquely determine the isomorphism type of the $\SL$-module $V$.  The previous two sentences together prove the uniqueness claim in the result. \\ \\
The proof of the result involving $\uq$ is the same as the above argument, except use Lemma \ref{thm:uq2mods} in place of Lemma \ref{thm:sl2mods}, and replace $V(d_j)$ with $V(d_j, 1)$.
\hfill $\Box$ \\

\noindent
We are now ready to prove Theorem \ref{thm:class}. \\

\noindent
{\it Proof of Theorem \ref{thm:class}:} \\
$(\Longrightarrow)$:  Adopt Assumption \ref{assume}.  We prove that the scalars in (\ref{parameterarray}) satisfy Theorem \ref{thm:class}(i)--(v).  \\ \\
First assume that $d=0$.  Theorem \ref{thm:class}(i),(ii),(v) hold, since they are vacuously true.  Theorem \ref{thm:class}(iv) holds, since it is trivially true.  Recall $V_0$ is an eigenspace of $A$, and by Definition \ref{def:shape}, $\rho_0 = \dim(V_0)$.  Thus, $\rho_0$ is a positive integer, and Theorem \ref{thm:class}(iii) holds. \\ \\
Now assume that $d \geq 1$.  By definition of eigenvalue (resp.~dual eigenvalue) sequence the sequence $\{ \theta_i \}_{i=0}^d$ (resp.~$\{ \theta^*_i \}_{i=0}^d$) is a list of the {\it distinct} eigenvalues of $A$ (resp.~$A^*$).  Thus, Theorem \ref{thm:class}(i) holds.  For $d=1$, Theorem \ref{thm:class}(ii) holds, since it is vacuously true.  For $d \geq 2$, Theorem \ref{thm:class}(ii) holds by Lemma \ref{thm:main2prep}.  Recall for $0 \leq i \leq d$, $V_i$ is an eigenspace of $A$, and by Definition \ref{def:shape}, $\rho_i = \dim(V_i)$.  Thus, $\rho_i$ is a positive integer for $0 \leq i \leq d$, and Theorem \ref{thm:class}(iii) holds.  By Definition \ref{def:shape}, $\rho_i = \dim(V_i) = \dim(V_{d-i}) = \rho_{d-i}$ for $0 \leq i \leq d$, and so Theorem \ref{thm:class}(iv) holds.  By Lemma \ref{thm:injsurj}, with $0 \leq i < d/2$ and $j=1$, the restriction $[A,A^*] |_{V_i} : V_i \rightarrow V_{i+1}$ is an injection.  From this and Definition \ref{def:shape}, $\rho_i = \dim(V_i) \leq \dim(V_{i+1}) = \rho_{i+1}$ for $0 \leq i < d/2$, and so Theorem \ref{thm:class}(v) holds. \\ \\
\noindent
$(\Longleftarrow)$:  Let $d$ denote a nonnegative integer, and let (\ref{parameterarray}) denote a sequence of scalars taken from $\K$ which satisfy Theorem \ref{thm:class}(i)--(v).  We construct a bidiagonal pair $A, \, A^*$ over $\K$ which has parameter array (\ref{parameterarray}).  For $d \leq 1$, define $b := 1$.  For $d \geq 2$, define $b$ as the common value of the expressions in Theorem \ref{thm:class}(ii).  We now break the argument into two cases depending of the value of $b$. \\ \\
Case 1:  $b=1$.  Let $V$ denote the segregated $\SL$-module constructed in Lemma \ref{thm:bijcorrsl2}.  For $0 \leq i \leq d$, let  $U_i$ denote the eigenspace for the action of $h$ on $V$ corresponding to the eigenvalue $d-2i$.  By Lemma \ref{thm:bijcorrsl2},
\begin{eqnarray}
\label{bijcorrsl21}
\rho_i = \dim(U_i) \qquad (0 \leq i \leq d).
\end{eqnarray}
Let $Y, \, Z$ denote an equitable pair in $\SL$.   By Theorem \ref{thm:uqsl2bidiag}, $Y, \, Z$ acts on $V$ as a reduced bidiagonal pair with base $1$.  We now show that
\begin{eqnarray}
\label{bijcorrsl22}
\mbox{the} \,\,  \mbox{parameter} \,\, \mbox{array} \,\,  \mbox{of} \,\, Y, \, Z \,\, \mbox{is} \,\, ( \{ 2i-d \}_{i=0}^d; \{ d-2i \}_{i=0}^d; \{ \rho_i \}_{i=0}^d).
\end{eqnarray}
The eigenvalue (resp.~dual eigenvalue) sequence of $Y, \, Z$ is $\{ 2i-d \}_{i=0}^d$ (resp.~$\{ d-2i \}_{i=0}^d$), since $Y, \, Z$ is reduced with base $1$.  We now show that the shape of $Y, \, Z$ is $\{ \rho_i \}_{i=0}^d$.  For $0 \leq i \leq d$, let  $Z_i$ denote the eigenspace for the action of $Z$ on $V$ corresponding to eigenvalue $d-2i$.  Identify the copy of $\SL$ given in Definition \ref{def:usl2} with the copy given in Theorem \ref{thm:ueq}, via the automorphism given in Theorem \ref{thm:ueq}.  Under this automorphism, $Z$ is mapped to $h$.  From this we have $\dim(Z_i) = \dim(U_i)$ for $0 \leq i \leq d$.  From this and (\ref{bijcorrsl21}), we have $\dim(Z_i) = \rho_i$ for $0 \leq i \leq d$.  So the shape of $Y, \, Z$ is $\{ \rho_i \}_{i=0}^d$.  We have now shown (\ref{bijcorrsl22}).  By Lemma \ref{thm:main2prep2}, there exists scalars $b_1, \, b_2, \, c_1, \, c_2$ in $\K$ with $b_2, \, c_2$ both nonzero such that
\begin{eqnarray}
\label{bijcorrsl23}
\theta_i = b_1 + b_2 2i, \qquad \theta^*_i = c_1 + c_2 (-2i)  \qquad (0 \leq i \leq d).
\end{eqnarray}
Define the linear transformations $A: V \rightarrow V$ and $A^* : V \rightarrow V$ as follows:  $A :=b_2 Y + (b_2 d + b_1) I$, $A^*:=c_2 Z + (c_1 -  c_2 d) I$.  Combining Lemma \ref{affineisbidiag}, (\ref{bijcorrsl22}), and (\ref{bijcorrsl23}), we find that $A, \, A^*$ is a bidiagonal pair on $V$ with parameter array $( \{ \theta_i \}_{i=0}^d; \{  \theta^*_i \}_{i=0}^d; \{ \rho_i \}_{i=0}^d)$.  \\ \\
Case 2:  $b \neq 1$.  Let $q$ denote a scalar in $\K$ such that $b=q^{-2}$.  Let $V$ denote the segregated $\uq$-module constructed in Lemma \ref{thm:bijcorrsl2}.  For $0 \leq i \leq d$, let  $U_i$ denote the eigenspace for the action of $k$ on $V$ corresponding to the eigenvalue $q^{d-2i}$.  By Lemma \ref{thm:bijcorrsl2},
\begin{eqnarray}
\label{bijcorrsl24}
\rho_i = \dim(U_i) \qquad (0 \leq i \leq d).
\end{eqnarray} 
Let $y, \, z$ denote an equitable pair in $\uq$.   By Theorem \ref{thm:uqsl2bidiag}, $y, \, z$ acts on $V$ as a reduced bidiagonal pair with base not equal to $1$.  We now show that
\begin{eqnarray}
\label{bijcorrsl25}
\mbox{the} \,\,  \mbox{parameter} \,\, \mbox{array} \,\,  \mbox{of} \,\, y, \, z \,\, \mbox{is} \,\, ( \{ q^{d-2i} \}_{i=0}^d; \{ q^{2i-d} \}_{i=0}^d; \{ \rho_i \}_{i=0}^d).
\end{eqnarray}
The eigenvalue (resp.~dual eigenvalue) sequence of $y, \, z$ is $\{ q^{d-2i} \}_{i=0}^d$ (resp.~$\{ q^{2i-d} \}_{i=0}^d$), since $y, \, z$ is reduced with base not equal to $1$.  We now show that the shape of $y, \, z$ is $\{ \rho_i \}_{i=0}^d$.  For $0 \leq i \leq d$, let  $y_i$ denote the eigenspace for the action of $y$ on $V$ corresponding to eigenvalue $q^{d-2i}$.  Identify the copy of $\uq$ given in Definition \ref{def:uqsl2} with the copy given in Theorem \ref{thm:uqeq}, via the automorphism given in Theorem \ref{thm:uqeq}.  Under this automorphism, $x$ is mapped to $k$.  Let $\Omega : V \rightarrow V$ denote the invertible linear operator from Lemma \ref{thm:omega}, and recall $\Omega ^{-1} x = y \Omega ^{-1}$.  Combining the previous two sentences, we have $\dim(y_i) = \dim(U_i)$ for $0 \leq i \leq d$.  From this and (\ref{bijcorrsl24}), we have $\dim(y_i) = \rho_i$ for $0 \leq i \leq d$.  So the shape of $y, \, z$ is $\{ \rho_i \}_{i=0}^d$.  We have now shown (\ref{bijcorrsl25}).  By Lemma \ref{thm:main2prep2}, there exists scalars $b_1', \, b_2 ', \, c_1 ', \, c_2 '$ in $\K$ with $b_2 ', \, c_2 '$ both nonzero such that 
\begin{eqnarray}
\label{bijcorrsl26}
\theta_i = b_1 ' + b_2 ' q^{-2i}, \qquad \theta^*_i = c_1 ' + c_2 ' q^{2i} \qquad (0 \leq i \leq d).
\end{eqnarray}
Define the linear transformations $A: V \rightarrow V$ and $A^* : V \rightarrow V$ as follows:  $A := q^{-d} \, b_2 ' \,  y + b_1 ' \, I$, $A^*:= q^d \, c_2 ' \, z + c_1 ' \, I$.  Combining Lemma \ref{affineisbidiag}, (\ref{bijcorrsl25}), and (\ref{bijcorrsl26}), we find that $A, \, A^*$ is a bidiagonal pair on $V$ with parameter array $( \{ \theta_i \}_{i=0}^d; \{  \theta^*_i \}_{i=0}^d; \{ \rho_i \}_{i=0}^d)$. \\

\noindent
The uniqueness claim in Theorem \ref{thm:class} follows immediately from Lemma \ref{thm:isompa}.
 \hfill $\Box$ \\

\section{Closing remarks}

\noindent 
In this section we make some suggestions for future research, and comment on some ongoing research which is related to this paper.  \\

\noindent
Let $A, \, A^*$ denote a bidiagonal pair as in Definition \ref{def:bidiag}.  Let $\lbrace W_i \rbrace _{i=0}^d$ denote the subspaces from Definition \ref{def:W}.  Let $B$ denote the linear transformation from Definition \ref{def:B1}, and recall that $\lbrace W_i \rbrace _{i=0}^d$ are the eigenspaces of $B$.  Recall by Lemma \ref{thm:AAbWij}(iii) that 
\begin{eqnarray}
\label{triple1}
A W_i \subseteq W_i + W_{i+1} \qquad (0 \leq i \leq d). 
\end{eqnarray}
Using (\ref{A1Ab1B13}) it can be shown that
\begin{eqnarray}
\label{triple2}
B V_i \subseteq V_{i-1} + V_{i} \qquad (0 \leq i \leq d).
\end{eqnarray}
Recall by Lemma \ref{thm:AAbWij}(iv) that 
\begin{eqnarray}
\label{triple3}
A^* W_i \subseteq W_{i-1} + W_{i} \qquad (0 \leq i \leq d). 
\end{eqnarray}
Using (\ref{A1Ab1B12}) it can be shown that
\begin{eqnarray}
\label{triple4}
B V^*_i \subseteq V^*_{i-1} + V^*_{i} \qquad (0 \leq i \leq d).
\end{eqnarray}
From (\ref{triple1}), (\ref{triple2}) (resp.~(\ref{triple3}), (\ref{triple4})) we see that $A, \, B$ (resp.~$A^*, \, B$) act on each others eigenspaces in a bidiagonal fashion.  Thus, for the triple $A, \, A^*, \, B$ of diagonalizable linear transformations, any two of them give a pair of bidiagonal actions on each others eigenspaces.  This suggests the idea of a {\it bidiagonal triple}.  We offer the following suggestion for a future research project.  Give a formal definition of bidiagonal triples, classify them, and explore the connections between bidiagonal pairs and bidiagonal triples.  \\

\noindent
Recall Theorem \ref{thm:main1}, which states that every bidiagonal pair satisfies the fundamental bidiagonal relation given in (\ref{main1.1}).  In terms of future research, what can be said about a pair of diagonalizable linear transformations which satisfy (\ref{main1.1})?  How close are such transformations from being a bidiagonal pair?  \\

\noindent
With reference to the definition of a bidiagonal pair, consider the following irreducibility condition: "there does not exist a subspace $W$ of $V$ such that $AW \subseteq W$, $A^* W \subseteq W$,  $W \neq 0$,  $W \neq V$."  What can be proven if condition (iii) of Definition \ref{def:bidiag} is replaced by this irreducibility condition?  This irreducibility condition is part of the definition of a tridiagonal pair.   \\

\noindent
Paul Terwilliger's forthcoming paper "Finite-dimensional irreducible $\uq$-modules from the equitable point of view" contains a result which extends Theorem \ref{thm:main3} of this paper.  This result gives a characterization of $\uq$ which shows how the equitable presentation of $\uq$ comes up naturally as the solution to a problem in linear algebra.  \\

\noindent
The concept of a Hessenberg pair generalizes both that of a bidiagonal pair and a tridiagonal pair.  Hessenberg pairs were introduced and characterized in \cite{Godjali09}.  For further information on Hessenberg pairs, see \cite{Godjali10, Godjali12}.  

\section{Appendix}

\noindent
In this appendix we prove Lemma \ref{thm:dequalsdelta}, namely that the number of eigenspaces of one transformation in a bidiagonal pair equals the number of eigenspaces of the other transformation. \\

\noindent
The following two lemmas will be used in the proof of Lemma \ref{thm:dequalsdelta}.

\begin{lemma}
\label{thm:polyinApolyinA*}
With reference to Definition \ref{def:bidiag}, let $\theta_i$ (resp.~$\theta_{i}^*$) denote the eigenvalue of $A$ (resp.~$A^*$) corresponding to $V_i$ (resp.~$V_{i}^*$) for $0 \leq i \leq d$ (resp.~$0 \leq i \leq \delta$).  Let $v \in V_0$ and $w \in V_0^*$.  Then the following {\rm(i)},{\rm(ii)} hold.
\begin{enumerate}
\item[\rm (i)] There exist $v_i \in V_i$ for $0 \leq i \leq d-1$ such that 
\begin{eqnarray*}
[A,A^*]^d v = (\theta_1 - \theta_0)(\theta_2 - \theta_1) \cdots (\theta_d - \theta_{d-1}) (A^{* \, d} v - A^{* \, d-1} v_0 - \cdots - A^* v_{d-2} - v_{d-1}).
\end{eqnarray*}
\item[\rm (ii)] There exist $v^*_i \in V^*_i$ for $0 \leq i \leq \delta -1$ such that 
\begin{eqnarray*}
[A,A^*]^{\delta} w = -(\theta^*_1 - \theta^*_0)(\theta^*_2 - \theta^*_1) \cdots (\theta^*_{\delta} - \theta^*_{\delta-1}) (A^{\delta} w - A^{\delta-1} v^*_0 - \cdots - A v^*_{\delta-2} - v^*_{\delta-1}).
\end{eqnarray*}
\end{enumerate}
\end{lemma}

\noindent
{\it Proof:} (i) The proof is by induction on $d$.  The result is trivially true for $d=0$.  Now assume that $d \geq 1$.  For $1 \leq i \leq d$, abbreviate $\Delta_i := (\theta_1 - \theta_0)(\theta_2 - \theta_1) \cdots (\theta_i - \theta_{i-1})$.  By induction, there exist $v_i \in V_i$ for $0 \leq i \leq d-2$ such that 
\begin{eqnarray}
\label{DeltaGamma1}
[A,A^*]^{d-1} v = \Delta_{d-1} \, \Gamma,
\end{eqnarray}    
where $\Gamma := A^{* \, d-1} v - A^{* \, d-2} v_0 - \cdots - A^* v_{d-3} - v_{d-2}$.  Combining Lemma \ref{thm:raise}(ii) and (\ref{DeltaGamma1}) yields $\Gamma \in V_{d-1}$.  Using this and (\ref{eq:b2}), we find that there exists $v_{d-1} \in V_{d-1}$ such that  
\begin{eqnarray}
\label{DeltaGamma2}
A^* \Gamma - v_{d-1} \in V_d.
\end{eqnarray}   
Now we have 
\begin{align*}
[A,A^*]^d v &= [A,A^*] \, [A,A^*]^{d-1} v \\
&= [A,A^*] \, \Delta_{d-1} \, \Gamma \qquad \hbox{by} \,\, (\ref{DeltaGamma1}) \\
&= \Delta _{d-1} \, (AA^* - A^*A) \, \Gamma \\
&= \Delta _{d-1} \, (A - \theta_{d-1} I) \, A^* \Gamma \qquad \hbox{since} \,\, \Gamma \in V_{d-1} \\
&= \Delta _{d-1} \, (A - \theta_{d-1} I) \, ((A^* \Gamma - v_{d-1}) + v_{d-1}) \\
&= \Delta _{d-1} \, (A - \theta_{d-1} I) \, (A^* \Gamma - v_{d-1}) \qquad \hbox{since} \,\, v_{d-1} \in V_{d-1} \\
&= \Delta _{d-1} \, (\theta_d - \theta_{d-1}) \, (A^* \Gamma - v_{d-1}) \qquad \hbox{by} \,\, (\ref{DeltaGamma2}) \\
&= \Delta _{d} \, (A^{* \, d} v - A^{* \, d-1} v_0 - \cdots - A^* v_{d-2} - v_{d-1}),
\end{align*}   
and we obtain the desired result. \\
(ii) Similar to (i).  
\hfill $\Box$ \\

\begin{lemma}
\label{thm:linindependent}
With reference to Definition \ref{def:bidiag}, the following {\rm(i)},{\rm(ii)} hold.
\begin{enumerate}
\item[\rm (i)] $I, \, A^*, \, A^{* \, 2}, \, \ldots, \, A^{* \, d}$ are linearly independent.
\item[\rm (ii)] $I, \, A, \, A^2, \, \ldots, \, A^{\delta}$ are linearly independent.
\end{enumerate}
\end{lemma}

\noindent
{\it Proof:} (i) Suppose, towards a contradiction, that  $I, \, A^*, \, A^{* \, 2}, \, \ldots, \, A^{* \, d}$ are linearly dependent.  Then there exist scalars $c_0, \, c_1, \, \ldots, \, c_d$, not all zero, such that $c_d A^{* \, d} + \cdots + c_1 A^{*} + c_0 I = 0$.  Define $r := \max\{ \, i \, |\, 0 \leq i \leq d\, , \, c_i \neq 0 \,\}$, and observe 
\begin{eqnarray}
\label{linearcombo1}
A^{* \, r} = -c_r^{-1} (c_{r-1} A^{* \, r-1} + \cdots + c_1 A^{*} + c_0 I).
\end{eqnarray} 
Using (\ref{eq:b2}) we have 
 \begin{eqnarray}
\label{linearcombo2}
A^{* \, j} V_i \subseteq V_i + V_{i+1} + \cdots + V_{i+j}
\end{eqnarray} 
for $0 \leq i \leq d$ and $0 \leq j \leq d-i$.  Let $0 \neq v \in V_0$.  From (\ref{linearcombo1}) we find that $A^{* \, d}$ is a linear combination of $A^{* \, r-1},  \, \cdots, \, A^{*}, \, I$.  Combining this and (\ref{linearcombo2}) yields $A^{* \, d} v \in V_0 + \cdots + V_{r-1} \subseteq V_0 + \cdots + V_{d-1}$.  This, Lemma \ref{thm:polyinApolyinA*}(i), and (\ref{linearcombo2}) give $[A,A^*]^d v \in V_0 + \cdots + V_{d-1}$.  However,  $[A,A^*]^d v \in V_d$ by Lemma \ref{thm:raise}(ii).  Combining the previous two sentences with the fact that $\sum_{i=0}^d V_i$ is a direct sum, we have
\begin{eqnarray}
\label{linearcombo3}
[A,A^*]^d v = 0.
\end{eqnarray}  
By (\ref{eq:b4}), $[A, A^*]^{d} |_{V_0} : V_0 \rightarrow V_{d}$ is a bijection.  So $[A,A^*]^d v \neq 0$ since $v \neq 0$.  This contradicts (\ref{linearcombo3}), and we obtain the desired result.   \\
(ii) Similar to (i).
\hfill $\Box$ \\

\noindent
We are now ready to prove Lemma \ref{thm:dequalsdelta}. \\

\noindent
{\it Proof of Lemma \ref{thm:dequalsdelta}:}  Since $A^*$ is diagonalizable, the number of eigenspaces of $A^*$ equals the degree of the minimal polynomial of $A^*$.  Hence, the degree of the minimal polynomial of $A^*$ is $\delta + 1$.  From Lemma \ref{thm:linindependent}(i), we see that the degree of the minimal polynomial of $A^*$ is greater than or equal to $d+1$.  Combining the previous two sentences, we have 
\begin{eqnarray}
\label{ddelta1}
\delta \geq d.
\end{eqnarray} 
Since $A$ is diagonalizable, the number of eigenspaces of $A$ equals the degree of the minimal polynomial of $A$.  Hence, the degree of the minimal polynomial of $A$ is $d + 1$.  From Lemma \ref{thm:linindependent}(ii), we see that the degree of the minimal polynomial of $A$ is greater than or equal to $\delta+1$.  Combining the previous two sentences, we have 
\begin{eqnarray}
\label{ddelta2}
d \geq \delta.
\end{eqnarray}  
Combining (\ref{ddelta1}) and (\ref{ddelta2}) yields $d = \delta$. 
\hfill $\Box$ \\

\section{Acknowledgments}
The author would like to thank his Ph.D. thesis advisor Paul Terwilliger for his guidance and many helpful suggestions.  Many of the ideas in this paper were developed while the author was a graduate student at the University of Wisconsin -- Madison.  The author would also like to thank the referee for a careful reading of the paper, and for his/her many helpful suggestions.    

\bibliographystyle{amsplain}
\bibliography{mybib}

\providecommand{\bysame}{\leavevmode\hbox to3em{\hrulefill}\thinspace}
\providecommand{\MR}{\relax\ifhmode\unskip\space\fi MR }
\providecommand{\MRhref}[2]{%
  \href{http://www.ams.org/mathscinet-getitem?mr=#1}{#2}
}
\providecommand{\href}[2]{#2}
\begin{thebibliography}{10}

\bibitem{Al-Najjar11}
H.~Al-Najjar, \emph{{L}eonard pairs associated with the equitable generators of
  the quantum algebra ${U}_q(\mathfrak{sl}_2)$}, Linear Multilinear Algebra
  \textbf{59} (2011), 1127--1142.

\bibitem{Al-Najjar04}
H.~Al-Najjar and B.~Curtin, \emph{A family of tridiagonal pairs}, Linear
  Algebra Appl. \textbf{390} (2004), 369--384.

\bibitem{Al-Najjar05}
\bysame, \emph{A family of tridiagonal pairs related to the quantum affine
  algebra ${U}_q(\widehat{\mathfrak{sl}}_2)$}, Electron. J. Linear Algebra
  \textbf{13} (2005), 1--9.

\bibitem{Al-Najjarinpress}
\bysame, \emph{A bilinear form for tridiagonal pairs of $q$-{S}erre type},
  Linear Algebra Appl. \textbf{428} (2008), no.~11--12, 2688--2698.

\bibitem{Al-Najjar10}
\bysame, \emph{Leonard pairs from the equitable basis of ${\mathfrak{sl}}_2$},
  Electron. J. Linear Algebra \textbf{20} (2010), 490--505.

\bibitem{Askey79}
R.~Askey and J.A. Wilson, \emph{A set of orthogonal polynomials that generalize
  the {R}acah coefficients or $6-j$ symbols}, SIAM J. Math Anal. \textbf{10}
  (1979), 1008--1016.

\bibitem{Bannai84}
E.~Bannai and T.~Ito, \emph{Algebraic combinatorics {I}: Association schemes},
  Benjamin/Cummings, London, 1984.

\bibitem{Baseilhac05}
P.~Baseilhac, \emph{Deformed {D}olan-{G}rady relations in quantum integrable
  models}, Nuclear Phys. B \textbf{709} (2005), 491--521.

\bibitem{Baseilhac052}
\bysame, \emph{An integrable structure related with tridiagonal algebras},
  Nuclear Phys. B \textbf{705} (2005), 605--619.

\bibitem{Baseilhac06}
\bysame, \emph{A family of tridiagonal pairs and related symmetric functions},
  J. Phys. A \textbf{39} (2006), no.~38, 11773--11791.

\bibitem{Benkart04}
G.~Benkart and P.~Terwilliger, \emph{Irreducible modules for the quantum affine
  algebra ${U}_q(\widehat{\mathfrak{sl}}_2)$ and its {B}orel subalgebra}, J.
  Algebra \textbf{282} (2004), 172--194.

\bibitem{Benkart07}
\bysame, \emph{The universal central extension of the three-point
  $\mathfrak{sl}_2$ loop algebra}, Proc. Amer. Math. Soc. \textbf{135} (2007),
  no.~6, 1659--1668.

\bibitem{Benkart10}
\bysame, \emph{The equitable basis for $\mathfrak{sl}_2$}, Math. Z.
  \textbf{268} (2011), 535--557.

\bibitem{Curtin071}
B.~Curtin, \emph{Modular {L}eonard triples}, Linear Algebra Appl. \textbf{424}
  (2007), no.~2--3, 510--539.

\bibitem{Curtin072}
\bysame, \emph{Spin {L}eonard pairs}, Ramanujan J. \textbf{13} (2007),
  no.~1--3, 319--332.

\bibitem{Elduque07}
A.~Elduque, \emph{The ${S}_4$-action on the tetrahedron algebra}, Proc. Roy.
  Soc. Edinburgh Sect. A \textbf{137} (2007), no.~6, 1227--1248.

\bibitem{Funk-Neubauer07}
D.~Funk-Neubauer, \emph{Raising/lowering maps and modules for the quantum
  affine algebra ${U}_q(\widehat {\mathfrak{sl}}_2)$}, Comm Algebra \textbf{35}
  (2007), no.~7, 2140--2159.

\bibitem{Funk-Neubauer09}
\bysame, \emph{Tridiagonal pairs and the $q$-tetrahedron algebra}, Linear
  Algebra Appl. \textbf{431} (2009), no.~5--7, 903--925.

\bibitem{Godjali09}
A.~Godjali, \emph{{H}essenberg pairs of linear transformations}, Linear Algebra
  Appl. \textbf{431} (2009), no.~9, 1579--1586.

\bibitem{Godjali10}
\bysame, \emph{Thin {H}essenberg pairs}, Linear Algebra Appl. \textbf{432}
  (2010), no.~12, 3231--3249.

\bibitem{Godjali12}
\bysame, \emph{Thin {H}essenberg pairs and double {V}andermonde matrices},
  Linear Algebra Appl. \textbf{436} (2012), no.~9, 3018--3060.

\bibitem{Hartwig05}
B.~Hartwig, \emph{Three mutually adjacent {L}eonard pairs}, Linear Algebra Appl
  \textbf{408} (2005), 19--39.

\bibitem{Hartwig07}
\bysame, \emph{The tetrahedron algebra and its finite dimensional irreducible
  modules}, Linear Algebra Appl. \textbf{422} (2007), no.~1, 219--235.

\bibitem{HarTer07}
B.~Hartwig and P.~Terwilliger, \emph{The tetrahedron algebra, the {O}nsager
  algebra, and the $\mathfrak{sl}_2$ loop algebra}, J. Algebra \textbf{308}
  (2007), no.~2, 840--863.

\bibitem{Humph72}
J.~Humphreys, \emph{Introduction to {L}ie algebras and representation theory},
  Springer-Verlag, New York, NY, 1972.

\bibitem{ItoTanTer01}
T.~Ito, K.~Tanabe, and P.~Terwilliger, \emph{Some algebra related to {P}- and
  {Q}-polynomial association schemes}, DIMACS Ser. Discrete Math. Theoret.
  Comput. Sci. (Providence, RI), vol.~56, American Mathematical Society, 2001,
  pp.~167--192.

\bibitem{ItoTer04}
T.~Ito and P.~Terwilliger, \emph{The shape of a tridiagonal pair}, J. Pure
  Appl. Algebra \textbf{188} (2004), 145--160.

\bibitem{ItoTer071}
\bysame, \emph{$q$-{I}nverting pairs of linear transformations and the
  $q$-tetrahedron algebra}, Linear Algebra Appl. \textbf{426} (2007), no.~2--3,
  516--532.

\bibitem{ItoTer072}
\bysame, \emph{The $q$-tetrahedron algebra and its finite dimensional
  irreducible modules}, Comm. Algebra \textbf{35} (2007), no.~11, 3415--3439.

\bibitem{ItoTer073}
\bysame, \emph{Tridiagonal pairs and the quantum affine algebra
  ${U}_q(\widehat{\mathfrak{sl}}_2)$}, Ramanujan J. \textbf{13 (1--3)} (2007),
  39--62.

\bibitem{ItoTer074}
\bysame, \emph{Tridiagonal pairs of {K}rawtchouk type}, Linear Algebra Appl.
  \textbf{427} (2007), no.~2--3, 218--233.

\bibitem{ItoTer075}
\bysame, \emph{Two non-nilpotent linear transformations that satisfy the cubic
  $q$-{S}erre relations}, J. Algebra Appl. \textbf{6} (2007), no.~3, 477--503.

\bibitem{ItoTerinpress4}
\bysame, \emph{Finite dimensional irreducible modules for the three-point
  $\mathfrak{sl}\sb 2$ loop algebra}, Comm. Algebra \textbf{36} (2008), no.~12,
  4557--4598.

\bibitem{ItoTerinpress1}
\bysame, \emph{Distance regular graphs and the $q$-tetrahedron algebra},
  European J. Combin. \textbf{30} (2009), no.~3, 682--697.

\bibitem{ItoTerinpress2}
\bysame, \emph{Distance regular graphs of $q$-{R}acah type and the
  $q$-tetrahedron algebra}, Michigan Math. J. \textbf{58} (2009), no.~1,
  241--254.

\bibitem{ItoTer09}
\bysame, \emph{Tridiagonal pairs of $q$-{R}acah type}, J. Algebra \textbf{322}
  (2009), no.~1, 68--93.

\bibitem{ItoTerinpress0}
\bysame, \emph{The augmented tridiagonal algebra}, Kyushu J. Math. \textbf{64}
  (2010), no.~1, 82--144.

\bibitem{ItoTer105}
\bysame, \emph{How to sharpen a tridiagonal pair}, J. Algebra Appl. \textbf{9}
  (2010), no.~4, 543--552.

\bibitem{ItoTerWang06}
T.~Ito, P.~Terwilliger, and C.W. Weng, \emph{The quantum algebra
  ${U}_q(\mathfrak{sl}_2)$ and its equitable presentation}, J. Algebra
  \textbf{298} (2006), 284--301.

\bibitem{Jantzen96}
J.C. Jantzen, \emph{Lectures on quantum groups}, American Mathematical Society,
  Providence, RI, 1996.

\bibitem{Koekoek98}
R.~Koekoek and R.F. Swarttouw, \emph{The {A}skey scheme of hypergeometric
  orthogonal polynomials and its $q$-analog},  (1998), available at
  http://aw.twi.tudelft.nl/~koekoek/research.html.

\bibitem{Koelink96}
H.T. Koelink, \emph{Askey-{W}ilson polynomials and the quantum $su(2)$ group:
  survey and applications}, Acta Appl. Math \textbf{44} (1996), 295--352.

\bibitem{Koelink00}
\bysame, \emph{$q$-{K}rawtchouk polynomials as spherical functions on the
  {H}ecke algebra of type {$B$}}, Trans. Amer. Math. Soc. \textbf{352} (2000),
  4789--4813.

\bibitem{Koornwinder93}
T.H. Koornwinder, \emph{Askey-{W}ilson polynomials as zonal spherical functions
  on the $su(2)$ quantum group}, SIAM J. Math. Anal. \textbf{24} (1993),
  795--813.

\bibitem{Leonard82}
D.~Leonard, \emph{Orthogonal polynomials, duality, and association schemes},
  SIAM J. Math. Anal. \textbf{13} (1982), 656--663.

\bibitem{Miki10}
K.~Miki, \emph{Finite dimensional modules for the $q$-tetrahedron algebra},
  Osaka J. Math \textbf{47} (2010), 559--589.

\bibitem{Nomura051}
K.~Nomura, \emph{A refinement of the split decomposition of a tridiagonal
  pair}, Linear Algebra Appl. \textbf{403} (2005), 1--23.

\bibitem{Nomura052}
\bysame, \emph{Tridiagonal pairs and the {A}skey-{W}ilson relations}, Linear
  Algebra Appl. \textbf{397} (2005), 99--106.

\bibitem{Nomura053}
\bysame, \emph{Tridiagonal pairs of height one}, Linear Algebra Appl.
  \textbf{403} (2005), 118--142.

\bibitem{Nomura061}
K.~Nomura and P.~Terwilliger, \emph{The determinant of ${A}{A}\sp *-{A}\sp
  *{A}$ for a {L}eonard pair ${A},{A}\sp *$}, Linear Algebra Appl. \textbf{416}
  (2006), no.~2--3, 880--889.

\bibitem{Nomura062}
\bysame, \emph{Matrix units associated with the split basis of a {L}eonard
  pair}, Linear Algebra Appl. \textbf{418} (2006), no.~2--3, 775--787.

\bibitem{Nomura063}
\bysame, \emph{Some trace formulae involving the split sequences of a {L}eonard
  pair}, Linear Algebra Appl. \textbf{413} (2006), 189--201.

\bibitem{Nomura071}
\bysame, \emph{Affine transformations of a {L}eonard pair}, Electron. J. Linear
  Algebra \textbf{16} (2007), 389--417.

\bibitem{Nomura072}
\bysame, \emph{Balanced {L}eonard pairs}, Linear Algebra Appl. \textbf{420}
  (2007), no.~1, 51--69.

\bibitem{Nomura073}
\bysame, \emph{Linear transformations that are tridiagonal with respect to both
  eigenbases of a {L}eonard pair}, Linear Algebra Appl. \textbf{420} (2007),
  no.~1, 198--207.

\bibitem{Nomura074}
\bysame, \emph{The split decomposition of a tridiagonal pair}, Linear Algebra
  Appl. \textbf{424} (2007), no.~2--3, 339--345.

\bibitem{Nomurainpress1}
\bysame, \emph{Sharp tridiagonal pairs}, Linear Algebra Appl. \textbf{429}
  (2008), no.~1, 79--99.

\bibitem{Nomurainpress2}
\bysame, \emph{The structure of a tridiagonal pair}, Linear Algebra Appl.
  \textbf{429} (2008), no.~7, 1647--1662.

\bibitem{Nomura081}
\bysame, \emph{The switching element for a {L}eonard pair}, Linear Algebra
  Appl. \textbf{428} (2008), no.~4, 1083--1108.

\bibitem{Nomurainpress3}
\bysame, \emph{Towards a classification of the tridiagonal pairs}, Linear
  Algebra Appl. \textbf{429} (2008), no.~2--3, 503--518.

\bibitem{Nomura09}
\bysame, \emph{Transition maps between the $24$ bases for a {L}eonard pair},
  Linear Algebra Appl. \textbf{431} (2009), no.~5--7, 571--593.

\bibitem{Nomura092}
\bysame, \emph{Tridiagonal pairs and the $\mu$-conjecture}, Linear Algebra
  Appl. \textbf{430} (2009), no.~1, 455--482.

\bibitem{Nomura102}
\bysame, \emph{On the shape of a tridiagonal pair}, Linear Algebra Appl.
  \textbf{432} (2010), no.~2--3, 615--636.

\bibitem{Nomura10}
\bysame, \emph{Tridiagonal pairs of $q$-{R}acah type and the $\mu$-conjecture},
  Linear Algebra Appl. \textbf{432} (2010), no.~12, 3201--3209.

\bibitem{Nomura11}
\bysame, \emph{Tridiagonal matrices with nonnegative entries}, Linear Algebra
  Appl. \textbf{434} (2011), no.~12, 2527--2538.

\bibitem{Noumi92}
M.~Noumi and K.~Mimachi, \emph{Askey-{W}ilson polynomials as spherical
  functions on the $su_q(2)$ quantum group}, Quantum Groups, Lecture Notes in
  Math. (Berlin), vol. 1510, Springer, 1992, pp.~98--103.

\bibitem{Rosengren99}
H.~Rosengren, \emph{Multivariable orthogonal polynomials as coupling
  coefficients for {L}ie and quantum algebra representations}, Centre for
  Mathematical Sciences, Lund University, Sweden, 1999, Ph.D. Thesis.

\bibitem{Rosengren07}
\bysame, \emph{An elementary approach to the $6j$-symbols (classical, quantum,
  rational, trigonometric, and elliptic)}, Ramanujan J. \textbf{13} (2007),
  131--166.

\bibitem{Ter90}
P.~Terwilliger, \emph{The incidence algebra of a uniform poset}, Math and its
  applications \textbf{20} (1990), 193--212.

\bibitem{Ter92}
\bysame, \emph{The subconstituent algebra of an association scheme {I}}, J.
  Algebraic Combin. \textbf{1} (1992), 363--388.

\bibitem{Ter011}
\bysame, \emph{Two linear transformations each tridiagonal with respect to an
  eigenbasis of the other}, Linear Algebra Appl. \textbf{330} (2001), 149--203.

\bibitem{Ter012}
\bysame, \emph{Two relations that generalize the $q$-{S}erre relations and the
  {D}olan-{G}rady relations}, Physics and Combinatorics 1999 (Nagoya) (River
  Edge, NJ), World Scientific Publishing, 2001, pp.~377--398.

\bibitem{Ter02}
\bysame, \emph{{L}eonard pairs from 24 points of view}, Rocky Mountain J. Math.
  \textbf{32} (2002), 827--888.

\bibitem{Ter03}
\bysame, \emph{Introduction to {L}eonard pairs}, J. Comput. Appl. Math.
  \textbf{153} (2003), no.~2, 463--475.

\bibitem{Ter041}
\bysame, \emph{{L}eonard pairs and the $q$-{R}acah polynomials}, Linear Algebra
  Appl. \textbf{387} (2004), 235--276.

\bibitem{Ter051}
\bysame, \emph{Two linear transformations each tridiagonal with respect to an
  eigenbasis of the other; comments on the parameter array}, Des. Codes
  Cryptogr. \textbf{34} (2005), 307--332.

\bibitem{Ter052}
\bysame, \emph{Two linear transformations each tridiagonal with respect to an
  eigenbasis of the other; comments on the split decomposition}, J. Comput.
  Appl. Math. \textbf{178} (2005), 437--452.

\bibitem{Ter053}
\bysame, \emph{Two linear transformations each tridiagonal with respect to an
  eigenbasis of the other; the ${T}{D}$-${D}$ and the ${L}{B}$-${U}{B}$
  canonical form}, J. Algebra \textbf{291} (2005), 1--45.

\bibitem{Ter042}
\bysame, \emph{An algebraic approach to the {A}skey scheme of orthogonal
  polynomials}, {O}rthogonal polynomials and special functions, Lecture Notes
  in Math. (Berlin), vol. 1883, Springer, 2006, pp.~255--330.

\bibitem{Terinpress}
\bysame, \emph{The equitable presentation for the quantum group
  ${U}_q({\mathfrak{g}})$ associated with a symmetrizable {K}ac-{M}oody algebra
  $\mathfrak{g}$}, J. Algebra \textbf{298} (2006), no.~1, 302--319.

\bibitem{TerVid04}
P.~Terwilliger and R.~Vidunas, \emph{{L}eonard pairs and the {A}skey-{W}ilson
  relations}, J. Algebra Appl. \textbf{3} (2004), 411--426.

\bibitem{Vidarinpress}
M.~Vidar, \emph{Tridiagonal pairs of shape (1,2,1)}, Linear Algebra Appl.
  \textbf{429} (2008), no.~1, 403--428.

\bibitem{Vid07}
R.~Vidunas, \emph{Normalized {L}eonard pairs and {A}skey-{W}ilson relations},
  Linear Algebra Appl. \textbf{422} (2007), no.~1, 39--57.

\bibitem{Vid08}
\bysame, \emph{{A}skey-{W}ilson relations and {L}eonard pairs}, Discrete Math.
  \textbf{308} (2008), no.~4, 479--495.

\bibitem{Zhedanov02}
A.S. Zhedanov and A.~Korovnichenko, \emph{Leonard pairs in classical
  mechanics}, J. Phys. A \textbf{5} (2002), 5767--5780.

\end{thebibliography}

\noindent
Darren Funk-Neubauer \hfil\break
\noindent Department of Mathematics and Physics \hfil\break
\noindent Colorado State University - Pueblo \hfil\break
\noindent 2200 Bonforte Boulevard \hfil\break
\noindent Pueblo, CO 81001 USA \hfil\break
email:   {\tt darren.funkneubauer@colostate-pueblo.edu} \hfil\break
phone:  (719) 549 - 2693 \hfil\break
fax:  (719) 549 - 2962

\end{document}